\documentclass[11pt]{article}
\usepackage{hyperref} 
\usepackage{marvosym}
\usepackage{graphicx}
\usepackage{latexsym,amsmath,amsfonts,amscd, amsthm, dsfont}
\usepackage{bm,color}
\usepackage{epsfig,verbatim,epstopdf,graphics}
\usepackage{subfigure}
\usepackage{changebar}
\usepackage{multirow}

\usepackage{algorithmic}

\usepackage{yhmath}
 \usepackage{booktabs} 
 \usepackage{tikz}
\usepackage{verbatim}
\usepackage{diagbox}
\usetikzlibrary{arrows,backgrounds,snakes,shapes}
 \numberwithin{equation}{section}

\graphicspath{{./}{./figure/}}
\allowdisplaybreaks

\newtheoremstyle{plainNoItalics}{}{}{\normalfont}{}{\bfseries}{.}{ }{}

\theoremstyle{plain}
\newtheorem{thm}{Theorem}[section]

\theoremstyle{plainNoItalics}

\newtheorem{rem}[thm]{Remark}
\newtheorem{prop}[thm]{Proposition}
\newtheorem{exa}[thm]{Example}

\newcommand{\bit}{\begin{itemize}}
\newcommand{\eit}{\end{itemize}}
\newcommand{\beq}{\begin{equation}}
\newcommand{\eeq}{\end{equation}}
\newcommand{\be}{\begin{eqnarray}}
\newcommand{\ee}{\end{eqnarray}}
\newcommand{\beno}{\begin{eqnarray*}}
\newcommand{\eeno}{\end{eqnarray*}}

\newcommand{\email}[1]{\protect\href{mailto:#1}{#1}}

\makeatletter

\newcommand{\Rmnum}[1]{\expandafter\@slowromancap\romannumeral #1@}

\makeatother

\begin{document}

\baselineskip=1.8pc


\title{A Generalized Eulerian-Lagrangian Discontinuous Galerkin Method for Transport Problems
 \thanks{
  Research of the first author is supported by the China Scholarship Council for 2 years' study at the University of Delaware.
Research of the second author is supported by NSF grant NSF-DMS-1818924, Air Force Office of Scientific Research FA9550-18-1-0257 and University of Delaware.
 }}
\author{Xue Hong
  \thanks{School of Mathematical Sciences, University of Science and Technology of China, Hefei, Anhui, 230026, P.R. China. (\email{xuehong1@mail.ustc.edu.cn}).}
  \and
  Jing-Mei Qiu%
  \thanks{Department of Mathematical Sciences, University of Delaware, Newark, DE, 19716, USA. (\email{jingqiu@udel.edu}).}
}
\maketitle
 \begin{abstract}
We propose a generalized Eulerian-Lagrangian (GEL) discontinuous Galerkin (DG) method. The method is a generalization of the Eulerian-Lagrangian (EL) DG method for transport problems proposed in [arXiv preprint arXiv: 2002.02930 (2020)],
 which tracks solution along approximations to characteristics in the DG framework, allowing extra large time stepping size with stability.
The newly proposed GEL DG method in this paper is motivated for solving linear hyperbolic systems with variable coefficients, where the velocity field for adjoint problems of the test functions is frozen to constant. In this paper, in a simplified scalar setting, we propose the GEL DG methodology by freezing the velocity field of adjoint problems, and by formulating the semi-discrete scheme over the space-time region partitioned by linear lines approximating characteristics. The fully-discrete schemes are obtained by method-of-lines Runge-Kutta methods. We further design flux limiters for the schemes to satisfy the discrete geometric conservation law (DGCL) and maximum principle preserving (MPP) properties. Numerical results on 1D and 2D linear transport problems are presented to demonstrate great properties of the GEL DG method. These include the high order spatial and temporal accuracy, stability with extra large time stepping size, and satisfaction of DGCL and MPP properties.
 \end{abstract}

\section{Introduction}

In this paper, we propose a generalized Eulerian-Lagrangian (GEL) Runge-Kutta (RK) discontinuous Galerkin (DG) method for a model transport equation in the form of
 \begin{equation}
u_t + \nabla \cdot( \mathbf{P}( u;\mathbf{x},t ) u ) = 0, \ (\mathbf{x},t)\in\mathbb{R}^d\times[0,T],
\label{general_nonlinear}
\end{equation}
where $d$ is the spatial dimension, $u: \mathbb{R}^d  \times[0,T]\rightarrow  \mathbb{R}$, and $\mathbf{P}( u;\mathbf{x},t  ) = ( P_1(u;\mathbf{x},t ), \cdots,P_d(u;\mathbf{x},t ) )^T$ is a linear or nonlinear velocity field. Such a model could come from a wide range of application fields including fluid dynamics, climate modeling, and kinetic description of plasma.

The GEL DG method is a generalization from the EL DG method proposed in \cite{cai2020eldg}. With RK time discretization, their fully-discrete schemes are termed ``GEL RK DG" and ``EL RK DG" methods.
The EL DG method is built upon a fixed set of computational mesh, yet in each time step, the solution is evolved over a local dynamic space-time region $\Omega_j$ (see Figure~\ref{iso1}), the partition of which is determined by linear approximations to characteristics. The EL DG method introduces modified adjoint problems for the test function
\beq
\label{eq: ad}
\phi_t +  \mathbf{\tilde{P}} \cdot \nabla \phi =0,
\eeq
with its velocity field $\mathbf{\tilde{P}}$ being a linear approximation to the velocity field $\mathbf{P}$ in \eqref{general_nonlinear}.  
Then RK methods are used for the time discretization via the method-of-lines approach. The proposed GEL DG method shares the same space-time partition strategy as the EL DG method. A major difference is the velocity field for the local modified adjoint problem of test functions. The GEL DG uses a constant function with $\mathbf{\tilde{P}} \equiv \bar{P}_j$ in \eqref{eq: ad}, whereas $\mathbf{\tilde{P}} \in P^1 (x, t)$ in the EL DG, to approximate the velocity field $\mathbf{P}$ of \eqref{general_nonlinear}. Such design is motivated from solving the wave equation via tracking information along different characteristics families in a linear system.
Take a 1-D 2-by-2 hyperbolic system for example,
\[
U_t + (A(x) U)_x = 0,
 \]
 which could arise from the wave propagation in a heterogeneous media. The modified adjoint problem for the system by the GEL DG method is
\[
\Phi_t + A_j \Phi_x = 0,
\]
where $A_j$ is a frozen local {\em constant} matrix that approximates $A(x)$ on $\Omega_j$.
In this paper, we focus on the GEL RK DG algorithm for scalar transport problems. The proposed GEL RK DG
maintains mass conservation, high order spatial and temporal accuracy, and allows for extra large time steps with stability. We also establish that the semi-discrete GEL DG and EL DG formulation are mathematically equivalent, whereas the time discretization introduces differences for fully discrete schemes. 

 We further study the property of discrete geometric conservation law (DGCL) and maximum principle preserving (MPP) property of GEL RK DG method, and find that the method fail to satisfy the DGCL and MPP property in general. We then propose MPP limiters to preserve the DGCL and MPP properties. The MPP limiters involve the polynomial rescaling limiter \cite{Zhang2010MPPlimiter,ZSpos} and the parametrized MPP flux limiter \cite{xiong2013parametrized,xiong2014high,Xiong2015MPPlimiter}.
The polynomial rescaling limiter preserves the MPP property for the piecewise DG polynomials, while the parametrized MPP flux limiter preserves the MPP property of cell averages in the final RK stage only, to avoid order reduction of RK methods if limiters are applied to intermediate RK solutions.

Finally, among different classes of EL methods in the literature, we would like to mention a few closely related ones. Eulerian-Lagrangian finite volume methods were introduced in \cite{huang2012eulerian} to handle nonlinearity for characteristic methods in the finite volume framework. The Eulerian Lagrangian Localized Adjoint Methods (ELLAM) \cite{celia1990eulerian} introduces an adjoint problem for the test function in the continuous finite element framework and has been applied to different problems \cite{wang1999family,russell2002overview}. Compared with ELLAM, the EL DG, SL DG and EL RK DG \cite{wang2007eulerian, cai2016high, cai2020eldg} are being developed in the discontinuous Galerkin finite element framework. Another line of development, that is closely related to this work, is the Arbitrary Lagrangian Eulerian (ALE) DG method \cite{klingenberg2017arbitrary,hong2020ALEDG_singular}. Both EL DG and ALE DG evolve the DG solution on a dynamic moving mesh. The dynamic mesh movement of the EL DG approximates characteristics for the potential of using larger time stepping sizes with stability, whereas the mesh movement of ALE DG could come from tracking moving computational domain and/or better shock resolution. The formulation of EL DG comes from the introduction of a local modified adjoint problem, whereas the ALE DG method is formulated through the coordinate transform of test function on a reference domain. 

This paper is organized as follows. In Section \ref{section:1d}, we develop the GEL DG for one-dimensional (1D) linear transport problems. We discuss numerical treatment of inflow boundary conditions and extensions to 2D problems by dimensional splitting.
In Section \ref{section:Stability}, we establish the equivalence of the GEL DG and EL DG method in semi-discrete form.
In Section \ref{section:GCL and MPP}, we study the DGCL and MPP properties of the GEL RK DG method and propose a MPP limiter to preserve DGCL and MPP for fully discrete schemes.
In Section \ref{section:numerical}, the performance of the proposed method is shown through extensive
numerical tests. Finally, concluding remarks are made in Section \ref{section:conclusion}.

\section{GEL DG formulation for linear transport problems} \label{section:1d}

We propose the GEL DG method, which differs from the EL DG method \cite{cai2020eldg} in the design of the modified adjoint problem for test functions. In the EL DG method, the adjoint problem is uniquely determined from the partition of the space-time region $\Omega_j$; while in the GEL DG method, the adjoint problem is independent from the partition of the space-time region. Such design of adjoint problems offers flexibility in handling hyperbolic systems when characteristic decomposition varies in space.

\subsection {1D linear transport problems} \label{subsection:1D linear transport problems}
We start from a 1D linear transport equation in the following form
\begin{equation}
u_t+(a(x, t)u)_x = 0, \quad x\in[x_a, x_b].
\label{scalar1d}
\end{equation}
For simplicity, we assume periodic boundary conditions, and the velocity field $a(x, t)$ is a continuous function of space and time. We perform a partition of the computational domain $x_a=x_{\frac12}< x_{\frac32}<\cdots< x_{N+\frac12} =x_b$.
Let $I_j=[x_{j-\frac12}, x_{j+\frac12} ]$ denote an element of length $\Delta x_j=x_{j+\frac12}-x_{j-\frac12}$ and define $\Delta x=\max_{j}\Delta x_j.$
We define the finite dimensional approximation space, $V_h^k = \{ v_h:  v_h|_{I_j} \in P^k(I_j) \}$, where $P^k(I_j)$ denotes the set of polynomials of degree at most $k$.
We let $t^n$ be the $n$-th time level and $\Delta t =t^{n+1}-t^n$ to be the time-stepping size.

The scheme formulation is summarized in four steps. We will first partition the space-time region $\Omega_j$'s, then introduce a new modified adjoint problem for the test function $\psi(x, t)$. We then formulate the semi-discrete GEL DG scheme. Finally, we apply the method-of-lines RK method for the time marching.

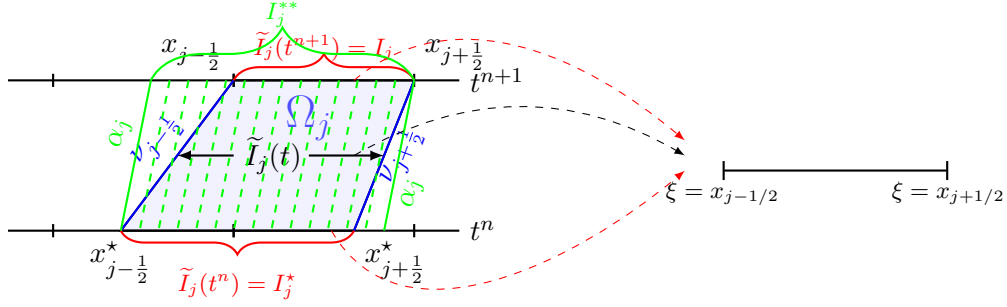
\begin{figure}[h!]
\label{iso1}
\centering
\begin{tikzpicture}[x=1cm,y=1cm]
  \begin{scope}[thick]

  \draw[fill=blue!5] (0.,2) -- (-1.5,0) -- (1.6,0) -- (2.4,2)
      -- cycle;
   \node[blue!70, rotate=0] (a) at ( 1. ,1.5) {\LARGE $\Omega_j$ };

   \draw (-3,3) node[fill=white] {};
    \draw (-3,-1) node[fill=white] {};
    \draw[black]                   (-3,0) node[left] {} -- (3,0)
                                        node[right]{$t^{n}$};
    \draw[black] (-3,2) node[left] {$$} -- (3,2)
                                        node[right]{$t^{n+1}$};

     \draw[snake=ticks,segment length=2.4cm] (-2.4,2) -- (0,2) node[above left] {$x_{j-\frac12}$};
     \draw[snake=ticks,segment length=2.4cm] (0,2) -- (2.4,2) node[above right] {$x_{j+\frac12}$};

          \draw[snake=ticks,segment length=2.4cm] (-2.4,0) -- (0,0);
     \draw[snake=ticks,segment length=2.4cm] (0,0) -- (2.4,0);

            \draw[blue,thick] (0.,2) node[left] {$$} -- (-1.5,0)
                                        node[black,below]{$x_{j-\frac12}^\star$ };

                \draw[blue,thick] (2.4,2) node[left] {$$} -- (1.6,0)
                                        node[black,below right]{ $x_{j+\frac12}^\star$};
\draw [decorate,color=red,decoration={brace,mirror,amplitude=9pt},xshift=0pt,yshift=0pt]
(-1.5,0) -- (1.6,0) node [red,midway,xshift=0cm,yshift=-20pt]
{\footnotesize $\widetilde{I}_j(t^{n}) = I_{j}^\star$};

\draw [decorate,color=red,decoration={brace,amplitude=10pt},xshift=0pt,yshift=0pt]
(0,2) -- (2.4,2) node [red,midway,xshift=0cm,yshift=13pt]
{\footnotesize $\widetilde{I}_j(t^{n+1}) = I_{j}$};

\draw [decorate,color=green,decoration={brace,amplitude=20pt},xshift=0pt,yshift=0pt]
(-1.1,2) -- (2.4,2) node [green,midway,xshift=-0.cm,yshift=25pt]
{\footnotesize $I^{**}_j$};
\draw[-latex,thick](0,1)node[right,scale=1.]{$\widetilde{I}_j(t)$}
        to[out=180,in=0] (-0.75,1) node[above left=2pt] {$$};  

\draw[-latex,thick](1,1)node[right,scale=1.]{ }
        to[out=0,in=180] (2,1) node[above left=2pt] {$$};  

  \node[blue, rotate=54] (a) at (-1.0,1.2) { $\nu_{j-\frac12}$ };
  \node[blue, rotate=70] (a) at (2.25,1.) { $\nu_{j+\frac12}$ };

\draw[green,thick] (-1.1,2) -- (-1.5,0);
\draw[green,thick] (2.4,2) -- (2.0,0);
\draw[green,dashed] (-0.85,2) -- (-1.25,0);
\draw[green,dashed] (-0.6,2) -- (-1.,0);
\draw[green,dashed] (-0.35,2) -- (-0.75,0);
\draw[green,dashed] (-0.1,2) -- (-0.5,0);
\draw[green,dashed] (0.15,2) -- (-0.25,0);
\draw[green,dashed] (0.4,2) -- (0.,0);
\draw[green,dashed] (0.65,2) -- (0.25,0);
\draw[green,dashed] (0.9,2) -- (0.5,0);
\draw[green,dashed] (1.15,2) -- (0.75,0);
\draw[green,dashed] (1.4,2) -- (1.,0);
\draw[green,dashed] (1.65,2) -- (1.25,0);
\draw[green,dashed] (1.9,2) -- (1.5,0);
\draw[green,dashed] (2.15,2) -- (1.75,0);

  \node[green, rotate=80] (a) at (-1.5,1.3) { $\alpha_j$ };
  \node[green, rotate=80] (a) at (2.38,0.5) { $\alpha_j$ };

  \end{scope}

   \draw[-latex,dashed](1.6,1)node[right,scale=1.0]{ }
        to[out=30,in=150] (6,1) node[] {};

     \draw[-latex,red,dashed](1.6,2)node[right,scale=1.0]{ }
        to[out=40,in=140] (6,1.2) node[] {};

     \draw[-latex,red,dashed](1.3,0)node[right,scale=1.0]{ }
        to[out=-60,in=220] (6,1-0.2) node[] {};

\draw[|-|,black,thick] (6.5,0.8) node[below]{\footnotesize $\xi=x_{j-1/2}$} to (9.5,0.8)node[below]{\footnotesize $\xi=x_{j+1/2}$};
\end{tikzpicture}
\caption{Illustration for dynamic element $\widetilde{I}_j(t)$ of new GEL DG. }
\end{figure}

\noindent
  {\bf (1) Partition of space-time region $\Omega_j$:} We define a space-time region $\Omega_j = \tilde{I}_j(t) \times [t^n, t^{n+1}]$ with
  $$\tilde{I}_j(t) = [\tilde{x}_{j-\frac12}(t), \tilde{x}_{j+\frac12}(t)], \quad t \in [t^n, t^{n+1}]$$
  being the dynamic interval, see Figure~\ref{iso1}. Here $\tilde{x}_{j\pm\frac12}(t) = x_{j\pm\frac12} + (t-t^{n+1}) \nu_{j\pm\frac12}$ are straight lines emanating from cell boundaries $x_{j\pm\frac12}$ with slopes $\nu_{j\pm\frac12}= a(x_{j\pm\frac12},t^{n+1})$ as in the EL DG scheme. We let ${I}^\star_j \doteq \tilde{I}_j(t^n) = [x^*_{j-\frac12}, x^*_{j+\frac12}]$ be the upstream cell of $I_j$ at $t^n$. Note that $\nu_{j\pm\frac12}$ are chosen to best take advantage of the characteristics information. In a classical Eulerian RK DG scheme, $\nu_{j+\frac12} =0$, $\forall j$.

\noindent {\bf (2) Adjoint problems.} We consider a local adjoint problem for the test function:
\begin{equation}
\begin{cases}
\psi_t + \alpha_j \psi_x =0 ,\ (x, t)\in\Omega_j,\\
\psi(t=t^{n+1}) = \Psi (x), \quad \forall \Psi(x) \in P^k(I^{**}_j).
\end{cases}
\label{eq: adjoint_new}
\end{equation}
Here we let $\alpha_j=a(x_j,t^{n+1})$.
To obtain the test function $\psi(x, t)$ on $\Omega_j$, $\Psi(x)$ needs to be defined in a large enough neighborhood containing $I_j$, named $I^{**}_j = [x^{**}_{j-\frac12},
x^{**}_{j+\frac12}]$, with
\beq
x^{**}_{j-\frac12} = \min(x_{j-\frac12}, x^*_{j-\frac12}+\alpha_j \Delta t), \qquad x^{**}_{j+\frac12} =\max(x_{j+\frac12}, x^*_{j+\frac12}+\alpha_j \Delta t).
\eeq
Please see the green curves in Figure~\ref{iso1} for the slope of $\alpha_j$, and the interval $I^{**}_j$. Here we take a natural extension of $\Psi$ from $I_j$ to $I^{**}_j$. The idea of defining $\Psi(x)$ on $I_j^{**}$ is to ensure that $\psi(x, t)$ can be properly found on $\Omega_j$, see the area shadowed by green dash lines in Figure~\ref{iso1}. This is different from the adjoint problem in the EL DG method, where the test function $\psi$ stays the same polynomial, if $\tilde{I}_j(t)$ is mapped to a reference interval $I_j$. 

\noindent{\bf (3) Formulation of the semi-discrete GEL DG scheme.} In order to formulate the scheme, we integrate $\text{\eqref{scalar1d}}\cdot\psi + \text{\eqref{eq: adjoint_new}}\cdot u$ over $\Omega_j$ ,
\begin{equation}
\int_{\Omega_j} \left[\text{\eqref{scalar1d}}\cdot\psi + \text{\eqref{eq: adjoint_new}}\cdot u \right] dxdt=0.
\end{equation}
That is,
\begin{align}
0 =& \int_{t^n}^{t^{n+1} } \int_{\tilde{I}_j(t)} \left( u_t\psi+u\psi_t \right)dxdt
+ \int_{t^n}^{t^{n+1} } \int_{\tilde{I}_j(t)}  \left( (a(x, t)u)_x\psi+\alpha_j \psi_x u \right) dxdt \nonumber \\
=&\int_{t^n}^{t^{n+1} } \int_{\tilde{I}_j(t)} ( u\psi )_t dxdt
 + \int_{t^n}^{ t^{n+1} } \int_{ \tilde{I}_j(t) } \left(  (a(x, t)u\psi)_x-a(x, t)u\psi_x + \alpha_j \psi_x u  \right)dxdt
 \nonumber \\
 =& \int_{t^n}^{t^{n+1} } \left[
\frac{d}{dt} \int_{\tilde{I}_j(t) }u\psi dx - \nu_{j+\frac12} u\psi |_{\tilde{x}_{j+\frac12}(t)}+\nu_{j-\frac12} u\psi |_{\tilde{x}_{j-\frac12}(t)} +   a u\psi\left|^{\tilde{x}_{j+\frac12}(t) }_{\tilde{x}_{j-\frac12}(t)} \right. \right.\nonumber\\
&\left.+\int_{ \tilde{I}_j(t) }  (\alpha_j-a)u \psi_x dx \right]dt
\nonumber\\
 =& \int_{t^n}^{t^{n+1} } \left[
\frac{d}{dt} \int_{\tilde{I}_j(t) }u\psi dx +(a_{j+\frac12}- \nu_{j+\frac12}) u\psi |_{\tilde{x}_{j+\frac12}(t)}-(a_{j-\frac12} -\nu_{j-\frac12}) u\psi |_{\tilde{x}_{j-\frac12}(t)} \right. \nonumber\\
& \left.
+\int_{ \tilde{I}_j(t) }  (\alpha_j-a)u \psi_x dx \right]dt.
\label{eq: int}
\end{align}
Letting $F(u) \doteq (a-\nu)  u$, the time differential form of \eqref{eq: int} gives
\begin{equation}
\frac{d}{dt} \int_{\tilde{I}_j(t)}(u\psi)dx =-  \left(F\psi\right) \left|_{\tilde{x}_{j+\frac12}(t) } \right. +   \left(F\psi\right) \left|_{\tilde{x}_{j-\frac12}(t) } \right.   + \int_{\tilde{I}_j(t)} (a-\alpha_j)u \psi_xdx.
\label{mol_2}
\end{equation}
Notice that the dynamic interval of $\tilde{I}_j(t)$ can always be linearly mapped to a reference cell {$\xi$ in $I_j$ by the mapping $\tilde{x}(t;(\xi,t^{n+1}))$ }.
then eq.~\eqref{mol_2} in the $\xi$-coordinate becomes
\begin{equation}
\frac{d}{dt} \int_{I_j}(u  \psi (\xi))\frac{\partial \tilde{x}(t;(\xi,t^{n+1})) }{\partial \xi}d\xi =-  \left(F\psi\right) \left|_{ \xi=x_{j+\frac12}  } \right. +   \left(F\psi\right) \left|_{ \xi=x_{j-\frac12}  } \right.      + \int_{I_j}(a-\alpha_j)u \psi_{\xi}d\xi.
\label{mol_4}
\end{equation}
The DG discretization of \eqref{mol_4} is to find $u_h(\xi, t) \in P^{k}(I_j)$,
so that the following equality holds,
\begin{equation}
\boxed{
\frac{d}{dt} \int_{I_j}u_h(\xi, t)  \psi (\xi, t)\frac{\partial \tilde{x} }{\partial \xi}d\xi =-  \hat{F}_{j+\frac12}\psi(\tilde{x}_{j+\frac12}^-(t), t) +   \hat{F}_{j-\frac12}\psi(\tilde{x}_{j-\frac12}^+(t), t)     + \int_{I_j}(a-\alpha_j)u_h \psi_\xi d\xi,
}
\label{1D GELDG scheme for scalar}
\end{equation}
for $\psi(x, t)$ satisfying the adjoint problem \eqref{eq: adjoint_new} with $\forall \Psi(x) = \psi(x, t^{n+1})\in P^k(I^{**}_j)$.
Here $\hat{F}$ at a cell boundary can be taken as a monotone flux, e.g. the Lax-Friedrichs flux
\beq
\hat{F}(u^-, u^+) = \frac12 (F(u^-) + F(u^+))+\frac{\alpha_0}{2}(u^- -u^+), \quad \alpha_0 = \max_{u} |F'(u)|;
\eeq
and we use Gauss quadrature rules with $k+1$ quadrature points to approximate the integral term on the R.S.H. of the equation \eqref{1D GELDG scheme for scalar}.

Next we discuss the choice of basis functions for representing solutions and test functions, with which one can assemble time-dependent mass matrices for implementation. In a classical Eulerian DG setting, the test functions are the same as basis functions of $V_h^k$. However, in the GEL DG setting, the test function $\psi_{j, m}(x, t)$ is in a time-dependent domain $\tilde{I}_j$ satisfying the adjoint problem ~\eqref{eq: adjoint_new} with
\beq
\label{eq: basis}
\Psi(x)=\Psi_{j,m}(x), \ j=1,...,N, \ m=0,...,k.
\eeq
$\{\Psi_{j,m}(x)\}_{1\leq j\leq N, 0\leq m\leq k}$ are the basis of $P^k(I_j)$ with a natural extension to $I_j^{**}$. In fact, we have from eq.~\eqref{eq: adjoint_new}
\beq
\label{eq: test function_t}
\psi_{j,m}(x,t) = \Psi_{j,m}(x-\alpha_j(t-t^{n+1})).
\eeq
We let $\tilde{U}_j(t)$ be a vector of size ${(k+1)\times 1}$ with its elements consisting of
 \begin{align}\label{tilde{u} for scalar}
   &\left\{ \int_{ \tilde{I}_j(t) } u_h(x,t) \psi_{j,m}(x,t) dx \right\}_{0 \leq m\leq k }.
 \end{align}
On the other hand, we set our basis function as $\left\{ \tilde{\psi}_{j,m}(x,t) \right\}_{1 \leq j\leq N,0 \leq m\leq k}$ on a reference cell with a mapping  $\frac{\tilde{x}(t;(\xi,t^{n+1}))-\tilde{x}_{j-\frac12}(t)}{\tilde{x}_{j+\frac12}(t)-\tilde{x}_{j-\frac12}(t)}=\frac{\xi-x_{j-\frac12}}{\Delta x_j}$ with
\beq
\label{eq: basis function_t}
\tilde{\psi}_{j,m}(\tilde{x}(t;(\xi,t^{n+1})),t) = \Psi_{j,m}(\xi),
\eeq
where $\Psi_{j,m}|_{I_j}$ is a set of basis in $V_h^k$ on $I_j$. Then we let
   \begin{align}\label{representation of basis}
    &u_h(x,t)=\sum_{l=0}^k \hat{u}^{(l)}_j(t)\tilde{\psi}_{j,l}(x,t), \quad \mbox{on} \quad \tilde{I}_j(t),
   \end{align}
   where $ \hat{u}^{(l)}$ are coefficients for the basis. Let $U_j(t) = (\hat{u}_j^{(0)}(t), \cdots, \hat{u}_j^{(k)}(t))^T$ be the coefficient vector of size $(k+1)\times 1$. Notice that $U_j(t)$ here is different from $\tilde{U}_j(t)$ defined in \eqref{tilde{u} for scalar}, satisfying
the GEL DG scheme \eqref{1D GELDG scheme for scalar}
   \begin{equation}
   \begin{aligned}\label{eqn of time dependent matrix element}
        &\int_{ \tilde{I}_j(t) } u_h(x,t) \psi_{j,m}(x,t) dx =\int_{ \tilde{I}_j(t) } \sum_{l=0}^k \hat{u}^{(l)}_j(t)\tilde{\psi}_{j,l}(x,t) \psi_{j,m}(x,t) dx\\
        &=\sum_{l=0}^k \hat{u}^{(l)}_j(t) \int_{ \tilde{I}_j(t) } \tilde{\psi}_{j,l}(x,t) \psi_{j,m}(x,t) dx=\sum_{l=0}^k \hat{u}^{(l)}_j(t) \int_{ \tilde{I}_j(t) } \tilde{\psi}_{j,l}(x,t) \Psi_{j,m}(x-\alpha_j (t-t^{n+1})) dx.
   \end{aligned}
   \end{equation}
   Now we assemble the time dependent mass-matrix $A_j(t)$ of size $(k+1)\times(k+1)$ with its elements consisting of
   $$\left\{ \int_{ \tilde{I}_j(t) } \tilde{\psi}_{j,l}(x,t) \Psi_{j,m}(x-\alpha_j (t-t^{n+1})) dx\right\}_{0 \leq l\leq k,0 \leq m\leq k}.$$
  Further, we have by eq.~\eqref{tilde{u} for scalar}, \eqref{representation of basis}, \eqref{eqn of time dependent matrix element}
   \beq\label{matrix representation of U}
   \tilde{U}_j(t)=A_j(t)\cdot U_j(t), \  \forall j=1,...,N, \ \forall t\in [t^n, t^{n+1}].
   \eeq
Now we can write the semi-discrete scheme \eqref{1D GELDG scheme for scalar} as
\begin{align}
&\frac{\partial}{\partial t} \tilde{U}_j(t)=\frac{\partial}{\partial t}(A_j(t)\cdot U_j(t)) = \mathcal{L}\left(U_{j-1}(t), U_j(t), U_{j+1}(t) ,t \right),
\label{semi}
\end{align}
where the spatial discretization operator on the RHS of \eqref{1D GELDG scheme for scalar} is denoted as $\mathcal{L}\left( U_{j-1}(t), U_j(t), U_{j+1}(t) ,t \right).$

\noindent{\bf (4) RK time discretization and fully discrete scheme.}\label{(4)fully discrete scheme}
Next, we describe the fully discrete scheme with method-of-lines RK discretization of the time derivative. There are two main steps involved here.
\begin{enumerate}
\item {\bf Obtain the initial condition} of \eqref{semi} by an $L^2$ projection of $u(x, t^n)$ from background cells onto upstream cells ${I}^*_j$. That is,
    \begin{align}
    &A_j(t^n)U_j^n=\left( \int_{ \tilde{I}_j(t^n) } u(x, t^n) \psi_{j,0}(x,t^n) dx, \cdots, \int_{ \tilde{I}_j(t^n) } u(x, t^n) \psi_{j,k}(x,t^n) dx  \right)^T\\
    &=\left( \int_{{I}_j^\star} u(x, t^n) \Psi_{j,0}(x+\alpha_j \Delta t^n) dx, \cdots, \int_{{I}_j^\star} u(x, t^n) \Psi_{j,0}(x+\alpha_j \Delta t^n) dx \right)^T.
    \end{align}
The integrals over the upstream cells above can be evaluated in the same fashion as the SL DG scheme \cite{cai2016high}.
%

\item {\bf Update \eqref{semi} from $U_j^n$ to $U_j^{n+1}$.} We apply the SSP explicit RK methods \cite{shu1988efficient} as in a method of lines approach.
In particular, the time-marching algorithm using an $s$-stage RK method follows the procedure below:
\begin{enumerate}
  \item Get the mesh information of the dynamic element $\tilde{I}_j^{(l)},l=0,\cdots, s$ on RK stages by $\tilde{x}_{j\pm\frac12}(t) = x_{j\pm\frac12} + (t-t^{n+1}) \nu_{j\pm\frac12}$.

  \item  For RK stages $i=1,\cdots, s$, let $t^{(l)}=t^n+d_l\Delta t^n$, compute
  \begin{align}
   A_j(t^{(i)})\cdot U_j^{(i)}= \sum_{l=0}^{i-1} \left[ \alpha_{il}  A_j(t^{(l)})\cdot U_j^{(l)}
  + \beta_{il} \Delta t \mathcal{L}\left( U_{j-1}^{(l)}, U_j^{(l)}, U_{j+1}^{(l)} ,t^{(l)} \right) \right],
  \end{align}

  where $\alpha_{il}$ and $\beta_{il}$ are related to RK methods, where  we can update the coefficients $U_j^{(i)}$ by inverting $A_j(t^{(i)})$ from the equation above. The coefficients for second, third and fourth order RK methods are provided in Table \ref{table:ssprk}.  
\end{enumerate}
\end{enumerate}
\begin{table}[!ht]
\caption{Parameters of some practical Runge-Kutta time discretizations. \cite{shu1988efficient,Jiang_Shu}}
\vspace{0.1in}
\centering
\begin{tabular}{ c  cc c }
\hline
Order  &   $\alpha_{il}$  & $\beta_{il}$  & $d_l$    \\
 \hline

 2  &    1  &  1 &  0\\
    & $\frac12$ \ $\frac12$ & 0 \ $\frac12$  &1 \\
\hline
 3  &    1  &  1 &  0\\
    & $\frac34$ \ $\frac14$ & 0 \ $\frac14$  &1 \\
    & $\frac13$ \ 0 \ $\frac23$  & 0 \ 0 \ $\frac23$ & $\frac12$ \\
\hline
 4  &    1                 &  $\frac12$               &  0\\
    & 1 \ 0                & 0 \ $\frac12$            &$\frac12$ \\
    & 1 \ 0 \ 0            & 0 \ 0 \ 1                & $\frac12$ \\
    & -$\frac13$ \ $\frac13$ \ $\frac23$ \ $\frac13$  & 0 \ 0 \ 0 \ $\frac16$       & 1 \\
\hline
\end{tabular}
\label{table:ssprk}
\end{table}
This finishes the description of a fully discrete GEL DG method, which enjoys the mass conservation as stated in the following Theorem.
\begin{thm} (Mass conservation)
Given a DG solution $u_h(x,t^n)\in V_h^k$ and assuming the boundary condition is periodic, the proposed fully discrete GEL DG scheme with SSP RK time discretization of \eqref{semi} is locally mass conservative.
In particular,
\begin{equation*}
\sum_{ i=1 }^N \int_{ I_j } u_h( x,t^{n+1} ) dx
=
\sum_{ i=1 }^N \int_{ I_j } u_h( x,t^n ) dx.
\end{equation*}
\end{thm}
\noindent
{\em Proof.} It can be proved by letting $\psi=1$, the conservative form of integrating $F$ function with unique flux at cell boundaries, as the mass conservation property of EL DG scheme \cite{cai2020eldg}. We skip details for brevity.

\subsection{Inflow Boundary conditions}
In this subsection, we discuss our treatment of inflow boundary conditions. We consider the linear transport equation \eqref{scalar1d} with the initial condition and the inflow boundary condition
\begin{equation}
\begin{cases}
u(x,0) = u_0(x),\\
u(x_b,t) = f(t).
\end{cases}
\label{eqn:1d}
\end{equation}
The proposed procedure for inflow boundary conditions follow steps below. At the outflow boundary, characteristics will go to the interior of domain, hence the original GEL DG algorithm could be directly applied.

\begin{figure}[h!]
\centering
\subfigure[]{
\begin{tikzpicture}[scale=0.5]
%

      \draw[fill=blue!10] (-2.5,3) -- (-6.,0) -- (-3.,0) -- (0,3)
      -- cycle;
      \draw[fill=green!10] (0,3) -- (-3.,0) -- (0.5,0) -- (2.5,3)
      -- cycle;
    \draw[black]                 (-7.5,0) node[left] { } -- (3,0)
                                        node[right]{\scriptsize$t^{n}$};
    \draw[black]                 (-7.5,3) node[left] { } -- (3,3)
                                        node[right]{\scriptsize$t^{n+1}$};

    \draw[black]                 (-2.5,0) node[left] { } -- (-2.5,3)
                                        node[right]{ };

%
%
%

    \draw[very thick] (0, 3-0.1) -- (0, 3+0.1) node[above] {\tiny$x_{ \frac32}$};
    \draw[very thick] (2.5, 3-0.1) -- (2.5, 3+0.1) node[above] {\tiny$x_{ \frac52}$};
    \draw[very thick] (-2.5, 3-0.1) -- (-2.5, 3+0.1) node[above left] {\tiny$x_{ b }$  };
    \draw[very thick] (-5, 3-0.1) -- (-5, 3+0.1) node[above] {\tiny$x_{ -\frac12}$ };
    \draw[very thick] (-7.5, 0-0.1) -- (-7.5, 0+0.1) node[above] { };
    \draw[very thick] (-7.5, 3-0.1) -- (-7.5, 3+0.1) node[above] {\tiny$x_{- \frac32}$ };
    \draw[very thick] (0, 0-0.1) -- (0, 0+0.1) node[above] { };
    \draw[very thick] (2.5, 0-0.1) -- (2.5, 0+0.1) node[above] { };
    \draw[very thick] (-2.5, 0-0.1) -- (-2.5, 0+0.1) node[above] { };
     \draw[ultra thick,red] (-6, 0-0.05) -- (-6, 0+0.05) node[red,below]{\tiny$x_{\frac12}^\star$ };
     \draw[ultra thick,red] (-3, 0-0.05) -- (-3, 0+0.05) node[red,below]{\tiny$x_{\frac32}^\star$ };
     \draw[ultra thick,red] (0.5, 0-0.05) -- (0.5, 0+0.05) node[red,below]{\tiny$x_{\frac52}^\star$ };

     \draw[ultra thick,black] (-2.5, 0-0.1) -- (-2.5, 0+0.1) node[above left ] {  };
     \draw[ultra thick,black] (-5, 0-0.1) -- (-5, 0+0.1) node[above left ] {  };

\node[blue!70, rotate=0] (a) at ( -3.0 ,1.5) {  $\Omega_{1}$ };

\node[green!100, rotate=0] (a) at ( 0 ,1.5) {  $\Omega_2$ };

\end{tikzpicture}
}
\subfigure[]{
\begin{tikzpicture}[scale=0.5]
\draw[white,fill=green!10] (9.5,2.3) to[out=200,in=80] (12-5,0)  -- (12-2.5,0)
      -- cycle;

\draw[white,fill=blue!10] (12-2.5,5) to[out=200,in=80] (12-7.5,0)  -- (12-5 ,0) to[out=70,in=220]  (12-2.5,2.3)
      -- cycle;
    \draw[black]                 (12-7.5,0) node[left] { } -- (12+3,0)
                                        node[right]{\scriptsize$t^{n}$};
    \draw[black]                 (12-7.5,3) node[left] { } -- (12+3,3)
                                        node[right]{\scriptsize$t^{n+1}$};

    \draw[black]                 (12-2.5,0) node[left] { } -- (12-2.5,5.2)
                                        node[right]{ };

\draw[-latex,dashed,red]( 12-5,0) node[right,scale=1.3]{$$}
        to[out=70,in=220] ( 12-2.5,2.3 ) node[right,scale=1.3] {\tiny$t_{ -\frac12 }^\star$ };
\draw[-latex,dashed,red]( 12-7.5,0) node[right,scale=1.3]{$$}
        to[out=80,in=200] ( 12-2.5,5 ) node[right,scale=1.3] {\tiny$t_{ -\frac32 }^\star$ };
    \draw[very thick] (12+0, 3-0.1) -- (12+0, 3+0.1) node[above] {\tiny$x_{ \frac32}$};
    \draw[very thick] (12+2.5, 3-0.1) -- (12+2.5, 3+0.1) node[above] {\tiny$x_{ \frac52}$};
    \draw[very thick] (12-2.5, 3-0.1) -- (12-2.5, 3+0.1) node[above left] {\tiny$x_{ b }$  };
    \draw[very thick] (12-5, 0-0.1) -- (12-5, 0.1) node[below] {\tiny$x_{ -\frac12}$};
    \draw[very thick] (12-7.5, 0-0.1) -- (12-7.5, 0.1) node[below] {\tiny$x_{ -\frac32}$};
    \draw[very thick] (12+0, 0-0.1) -- (12+0, 0+0.1) node[above] { };
    \draw[very thick] (12+2.5, 0-0.1) -- (12+2.5, 0+0.1) node[above] { };
    \draw[very thick] (12-2.5, 0-0.1) -- (12-2.5, 0+0.1) node[above] { };
 \draw[red]                 (12-5,0) node[left] { } -- (12-2.5,0)
                                        node[right]{ };
    \draw[ultra thick,red] (12-5, 0-0.1) -- (12-5, 0+0.1);

     \draw[ultra thick,red] (12-2.5, 0-0.1) -- (12-2.5, 0+0.1) node[above left ] {  };

     \fill [blue] ( 12-2.5, 2 ) circle (2pt) node[right] {};

     \fill [blue] ( 12-2.5, 1.5  ) circle (2pt) node[right] {};

     \fill [blue] ( 12-2.5, 1.5-1.1618955 ) circle (2pt) node[right] {};


     \fill [red] ( 12-3, 0 ) circle (1.8pt) node[right] {};
     \fill [red] ( 12-4, 0 ) circle (1.8pt) node[right] {};

      \fill [red] ( 12-4.5, 0 ) circle (1.8pt) node[right] {};

  \draw[-latex, blue]( 12-2.5, 2. )node[left,scale=1.3]{$$}
        to[out=220,in=75] ( 12-4.5,0) node[below=2pt] { };

    \draw[-latex, blue]( 12-2.5, 1.5 )node[left,scale=1.3]{$$}
        to[out=200,in=80] ( 12-4 ,0) node[below=2pt] { };

      \draw[-latex, blue]( 12-2.5, 1.5-1.1618955 )node[left,scale=1.3]{$$}
        to[out=190,in=70] ( 12-3 ,0) node[below=2pt] { };

        \node[blue!70, rotate=0] (a) at ( 12-6 ,1.5) {  $\Omega_{-1}^*$ };
        \node[blue!70, rotate=0] (a) at ( 12-3. ,1) {  $\Omega_{0}^*$ };

      \node[red!90, rotate=0] (a) at ( 12-3.7 ,-0.3) {  $I_b$ };

\end{tikzpicture}
}
\caption{Illustration on ghost cells intersecting inflow boundary.}
\label{boundaryproblem_1d}
\end{figure}
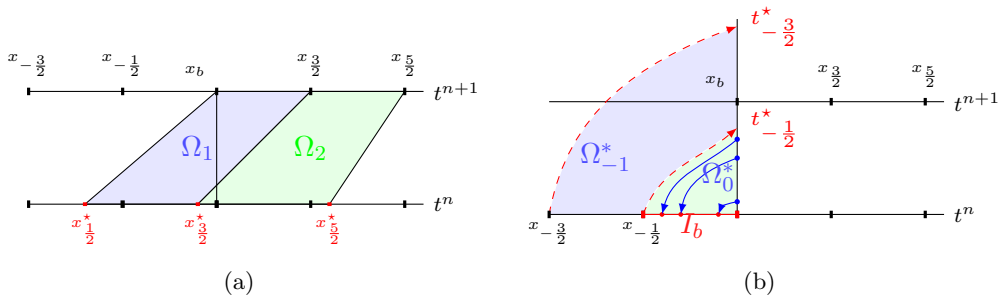

Step 1:
Set up ghost cells.
We first set up a ghost region which is sufficiently large, on which discretizations are performed to define ghost cells. For example, in the 1D setting, for $CFL<2$, we have two ghost cells $[x_{-\frac32},x_{-\frac12}]$, $[x_{-\frac12},x_{\frac12}]$ as illustrated in Figure \ref{boundaryproblem_1d} (a).
%

Step 2: Obtain the DG solutions on ghost cells.
We find DG solution at $t^n$ on ghost cells by tracking information along characteristics from boundary data in the semi-Lagrangian fashion. In order to do this, we first find the velocity field $a(x, t)$ on ghost regions, which could be done by a natural extrapolation from interior of domain. 
In particular, we consider the following problem
\begin{equation}
\begin{cases}
u_t + ( a(x,t)u)_x = 0, \quad \mbox{x at the ghost region}\\
u(x_b,t) = f(t),
\end{cases}
\label{eq: inflow_u}
\end{equation}
where $a(x, t)$ at the ghost region can be approximated by extrapolations from the interior of domain.
As shown in Figure \ref{boundaryproblem_1d}(b), there is $\Omega_{0}^*$ bounded by characteristic curves emanating from boundaries of ghost cells.
We let the test function $\psi(x,t)$ satisfies the adjoint problem with $\forall \Psi\in P^k(I_b)$,
\begin{equation}
\begin{cases}
\psi_t + a(x,t) \psi_x =0 ,\\
\psi(x, t^{n}) = \Psi(x), \ x\in[x_{-\frac12}, x_{\frac12}].
\end{cases}
\label{final-value}
\end{equation}
Integrate $\left(\eqref{eq: inflow_u} \cdot \psi + \eqref{final-value} \cdot u \right)$ over $\Omega_0^*$, we have
\begin{equation}
  \int_{\Omega_{0}^*} (u\psi)_t+(a(x,t)u\psi)_xdxdt=0.
\end{equation}
Using Green formula, we can get
\begin{equation}
\int_{x_{-\frac12}}^{x_{\frac12} } u(x,t^n)\Psi(x) dx
= \int_{t^n}^{t^*_{-\frac12}} a(x_b,t) u(x_b,t) \psi(x_b,t)dt.
\end{equation}
As in the SL DG \cite{cai2019comparison,cai2018high,cai2016high}, Gauss quadrature rule can be applied to evaluate the right-hand side of the above equation. Similar procedure can be used to obtain DG solutions on the ghost cells.

Step 3: Update solution. Once the solution on ghost cells are available, we can update the solution following the GEL RK DG procedure described previously.
%

\subsection {2D linear transport problems}

We extend the GEL RK DG algorithm to 2D problems via dimensional splitting \cite{qiu2011positivity}. Consider a linear 2D transport equation
\begin{equation}\label{eqn: 2d linear equation}
  u_t + (a(x,y,t)u)_x + (b(x,y,t)u)_y = 0, (x,y) \in \Omega,
\end{equation}
with a proper initial condition $u(x,y,0) = u_0(x,y)$ and boundary conditions. Here $(a(x,y,t), b(x,y,t))$ is a velocity field. The domain $\Omega$ is partitioned into rectangular meshes with each computational cell $A_{ij} = [x_{i- \frac{1}{2}}, x_{i+\frac{1}{2}}] \times [y_{j+ \frac{1}{2}}, y_{j+\frac{1}{2}}]$, where we use the piecewise $Q^k$ tensor-product polynomial spaces.
\begin{figure}[h!]
\centering
\begin{tikzpicture}[x=1cm,y=1cm]
  \begin{scope}[thick]

  \draw[fill=white!5] (-2.,-2.) -- (-2.,2) -- (2.,2.) -- (2.,-2.)
      -- cycle;
  \fill[blue] (-1.722,-1.722) circle (0.1);
  \fill[blue] (-1.722,-0.68) circle (0.1);
  \fill[blue] (-1.722,1.722) circle (0.1);
  \fill[blue] (-1.722,0.68) circle (0.1);
  \fill[blue] (-0.68,-1.722) circle (0.1);
  \fill[blue] (-0.68,-0.68) circle (0.1);
  \fill[blue] (-0.68,0.68) circle (0.1);
  \fill[blue] (-0.68,1.722) circle (0.1);
  \fill[blue] (0.68,-1.722) circle (0.1);
  \fill[blue] (0.68,-0.68) circle (0.1);
  \fill[blue] (0.68,0.68) circle (0.1);
  \fill[blue] (0.68,1.722) circle (0.1);
  \fill[blue] (1.722,-1.722) circle (0.1);
  \fill[blue] (1.722,-0.68) circle (0.1);
  \fill[blue] (1.722,1.722) circle (0.1);
  \fill[blue] (1.722,0.68) circle (0.1);
  \end{scope}

\begin{scope}[thick]
 \draw[fill=white!5] (3.,-2.) -- (3.,2) -- (7.,2.) -- (7.,-2.)
      -- cycle;
  \fill[blue] (3.278,-1.722) circle (0.1);
  \fill[blue] (3.278,-0.68) circle (0.1);
  \fill[blue] (3.278,1.722) circle (0.1);
  \fill[blue] (3.278,0.68) circle (0.1);
  \fill[blue] (4.32,-1.722) circle (0.1);
  \fill[blue] (4.32,-0.68) circle (0.1);
  \fill[blue] (4.32,0.68) circle (0.1);
  \fill[blue] (4.32,1.722) circle (0.1);
  \fill[blue] (5.68,-1.722) circle (0.1);
  \fill[blue] (5.68,-0.68) circle (0.1);
  \fill[blue] (5.68,0.68) circle (0.1);
  \fill[blue] (5.68,1.722) circle (0.1);
  \fill[blue] (6.722,-1.722) circle (0.1);
  \fill[blue] (6.722,-0.68) circle (0.1);
  \fill[blue] (6.722,1.722) circle (0.1);
  \fill[blue] (6.722,0.68) circle (0.1);
  \draw[red] (3.,-1.722) -- (7.,-1.722);
  \draw[red] (3.,-0.68) -- (7.,-0.68);
  \draw[red] (3.,0.68) -- (7.,0.68);
  \draw[red] (3.,1.722) -- (7.,1.722);
  \draw[red][->] (7.2,-1.722) -- (7.5,-1.722);
  \draw[red][->] (7.2,-0.68) -- (7.7,-0.68);
  \draw[red][->] (7.2,0.68) -- (7.9,0.68);
  \draw[red][->] (7.2,1.722) -- (8.1,1.722);
\end{scope}
\begin{scope}[thick]
 \draw[fill=white!5] (8.5,-2.) -- (8.5,2) -- (12.5,2.) -- (12.5,-2.)
      -- cycle;
  \fill[blue] (8.778,-1.722) circle (0.1);
  \fill[blue] (8.778,-0.68) circle (0.1);
  \fill[blue] (8.778,1.722) circle (0.1);
  \fill[blue] (8.778,0.68) circle (0.1);
  \fill[blue] (9.82,-1.722) circle (0.1);
  \fill[blue] (9.82,-0.68) circle (0.1);
  \fill[blue] (9.82,0.68) circle (0.1);
  \fill[blue] (9.82,1.722) circle (0.1);
  \fill[blue] (11.18,-1.722) circle (0.1);
  \fill[blue] (11.18,-0.68) circle (0.1);
  \fill[blue] (11.18,0.68) circle (0.1);
  \fill[blue] (11.18,1.722) circle (0.1);
  \fill[blue] (12.222,-1.722) circle (0.1);
  \fill[blue] (12.222,-0.68) circle (0.1);
  \fill[blue] (12.222,1.722) circle (0.1);
  \fill[blue] (12.222,0.68) circle (0.1);
  \draw[red] (8.778,-2.) -- (8.778,2.);
  \draw[red] (9.82,-2.) -- (9.82,2.);
  \draw[red] (11.18,-2.) -- (11.18,2.);
  \draw[red] (12.222,-2.) -- (12.222,2.);
  \draw[red][->] (8.778,2.2) -- (8.778,2.5);
  \draw[red][->] (9.82,2.2) -- (9.82,2.7);
  \draw[red][->] (11.18,2.2) -- (11.18,2.9);
  \draw[red][->] (12.222,2.2) -- (12.222,3.1);

\end{scope}
\end{tikzpicture}
\label{splittingk3}
\caption{Illustration of the 2D GEL RK DG scheme via Strang splitting. k = 3. }
\end{figure}
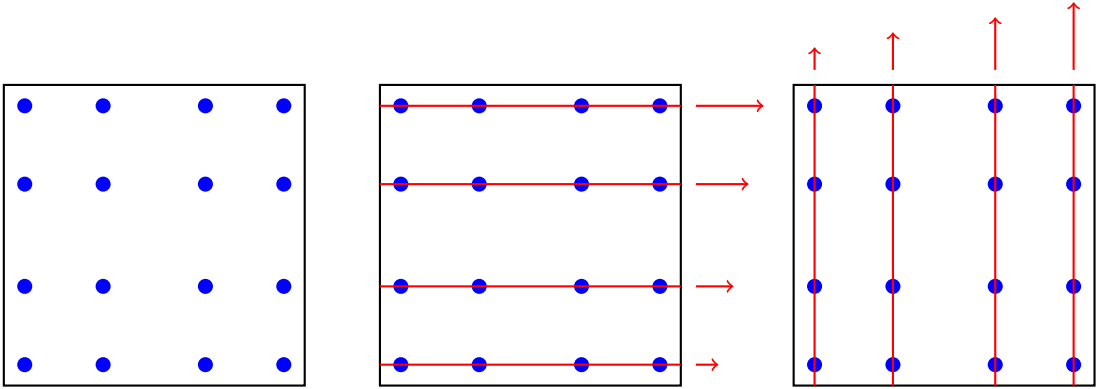
\begin{enumerate}
 \item We first locate $(k + 1)^2$ tensor-product Gaussian nodes on cell $A_{ij}: (x_{i,p}, y_{j,q}),\  p,q = 0,... ,k$. For example, see Figure \ref{splittingk3} (left) for the case of $k=3$.
 \item Then, the equation \eqref{eqn: 2d linear equation} is split into two 1D advection problems based on the quadrature nodes in both $x-$ and $y-$ directions:
\begin{align}\label{eqn: split the 2d linear equation 1}
  u_t + (a(x,y,t)u)_x &= 0,\\
 \label{eqn: split the 2d linear equation 2} u_t + (b(x,y,t)u)_y &= 0.
\end{align}
Based on a 1D GEL RK DG formulation, the split equations \eqref{eqn: split the 2d linear equation 1} and \eqref{eqn: split the 2d linear equation 2} are evolved via Strang splitting
over a time step $\Delta t$ as follows.
\begin{itemize}
  \item Evolve 1D equation \eqref{eqn: split the 2d linear equation 1} at different $y_{j,q}'s$ with different velocity for a half
time-step $\Delta t/2$, see Figure \ref{splittingk3} (middle). For each $y_{j,q}$, the $(k + 1)$ point values are mapped to a $P^k$ polynomial per cell, then the 1D equation \eqref{eqn: split the 2d linear equation 1} is evolved by the proposed GEL RK DG scheme. Finally, we can map the evolved $P^k$ polynomial back to the $(k + 1)$ point values to update the solution.
  \item Evolve 1D equation \eqref{eqn: split the 2d linear equation 2} at different $x_{i,p}'s$ for a full time-step $\Delta t$ as above, see Figure \ref{splittingk3} (right).
  \item Evolve 1D equation \eqref{eqn: split the 2d linear equation 1} at different $y_{j,q}'s$ for another half time-step $\Delta t/2$.
\end{itemize}
\end{enumerate}
  The splitting 2D GEL DG formulation maintains many desired properties from the base 1D formulation, such as high order accuracy in space, extra large time stepping size with stability and mass conservation. However, a second splitting error in time is introduced, and the computational cost increases exponentially with the dimension of the problem. Higher order splitting methods can be constructed in the spirit of composition methods \cite{yoshida,hairerGeom,Forest,cai2016conservative}. A direct 2D algorithm as those in \cite{cai2020eldg} will be subject to our future work.

\section{Equivalence of semi-discrete GEL DG and EL DG methods}\label{section:Stability}

In this section, we first study the equivalence between the evolution step of the GEL DG and SL DG methods for a linear constant problem \eqref{scalar1d} in section \ref{subsection: equivalence of GELDG and SLDG}; the equivalence between the GEL DG and EL DG methods for a linear variable coefficient problem \eqref{scalar1d} is also presented in section \ref{section:Equivalence between GELDG and ELDG}. In \cite{cai2020eldg}, we established that the evolution step in the EL DG scheme is the same as the ALE DG method, for which theoretical stability analysis are performed in \cite{klingenberg2017arbitrary}.
%
\subsection{Equivalence between semi-discrete GEL DG and SL DG for a linear constant coefficient equation}\label{subsection: equivalence of GELDG and SLDG}
The extra degree of freedom in the GEL DG scheme design, compared with the SL DG, is the space-time partition and the adjoint problem for the test function. For a linear constant coefficient advection equation, if we assume the exact space-time partition as the SL DG method, while varying the velocity field of the modified adjoint problem in GEL DG, we show below that the semi-discrete GEL DG scheme is equivalent to the SL DG scheme.
\begin{thm}\label{lemma:equivalence between GELDG and SLDG}
  For linear constant coefficient equation \eqref{scalar1d} with $a(x,t)=1$, GEL DG scheme with the exact space-time partition $\nu_{j\pm\frac12}=1$ and a perturbation of velocity $\alpha_j=1+c, c\neq 0$ in \eqref{eq: adjoint_new}, is equivalent to SL DG scheme in semi-discrete form.
  \begin{proof}
    In this case, the GEL DG scheme \eqref{mol_2} is reduced into
\begin{equation}\label{the semi-discrete scheme for linear scaler equation}
  \frac{d}{dt} \int_{\tilde{I}_j(t)}(u^{GEL DG}(x,t)\psi^{GEL DG}(x,t))dx = -\int_{\tilde{I}_j(t)} cu^{GEL DG} \psi^{GEL DG}_xdx.
\end{equation}
We can rewrite the scheme as in integral form
\begin{equation}\label{eqn:GELDG in integral form}
\begin{aligned}
\int_{{I}_j} &(u^{GELDG}(x,t^{n+1})\Psi(x))dx\\
&=\int_{\tilde{I}_j(t^n)}(u^{GELDG}(x,t^{n})\psi^{GELDG}(x,t^n))dx-c\int^{t^{n+1}}_{t^n}\int_{\tilde{I}_j(\tau)} u^{GELDG}(x,\tau) \psi^{GELDG}_x(x,\tau)dxd\tau\\
&=\int_{\tilde{I}_j(t^n)}(u^{n}(x)\Psi(x+(1+c)\Delta t))dx-c\int^{t^{n+1}}_{t^n}\int_{\tilde{I}_j(\tau)} u^{GELDG}(x,\tau) \Psi_x(x+(1+c)(t^{n+1}-\tau))dxd\tau\\
&\doteq RHS^{GELDG},
\end{aligned}
\end{equation}
where $\tilde{I}_j(t^n)=[x_{j-\frac12} -\Delta t,x_{j+\frac12} -\Delta t]$ and $\tilde{I}_j(t)=[x_{j-\frac12} + t-t^{n+1}, x_{j+\frac12} + t-t^{n+1}]$ .
We also rewrite the SL DG scheme as in integral form
\begin{equation}\label{eqn:SLDG in integral form}
  \begin{aligned}
  &\int_{{I}_j} (u^{SLDG}(x,t^{n+1})\Psi(x))dx=\int_{\tilde{I}_j(t^n)}(u^{n}(x)\psi^{SLDG}(x,t^n))dx\\
  &=\int_{\tilde{I}_j(t^n)}(u^{n}(x)\Psi(x+\Delta t))dx\doteq RHS^{SLDG},
  \end{aligned}
\end{equation}
If we assume $u^{GELDG}(x,\tau)=u^{SLDG}(x,\tau)$, $\tau\in[t^n, t^{n+1})$ then
\begin{align*}
  &RHS^{SLDG}-RHS^{GELDG}\\
  &=\int_{\tilde{I}_j(t^n)}(u^n(x)\Psi(x+\Delta t))dx-\int_{\tilde{I}_j(t^n)}(u^n(x)(x,t^{n})\Psi(x+(1+c)\Delta t))dx\\
  &+c\int^{t^{n+1}}_{t^n}\int_{\tilde{I}_j(t^n)} u^{SLDG}(\xi+(\tau-t^n),\tau) \Psi_x(\xi+c(t^{n+1}-\tau)+\Delta t)d\xi d\tau\\
  &=\int_{\tilde{I}_j(t^n)}u^n(x)(\Psi(x+\Delta t)-\Psi(x+(1+c)\Delta t))dx-\int_{\tilde{I}_j(t^n)}u^n(\xi) \int^{t^{n+1}}_{t^n}-c\Psi_x(\xi+c(t^{n+1}-\tau)+\Delta t)d\tau d\xi\\
  &=\int_{\tilde{I}_j(t^n)}u^n(x)(\Psi(x+\Delta t)-\Psi(x+(1+c)\Delta t))dx-\int_{\tilde{I}_j(t^n)}u^n(\xi) \int^{t^{n+1}}_{t^n}\Psi_{\tau}(\xi+c(t^{n+1}-\tau)+\Delta t)d\tau d\xi\\
  &=0,
\end{align*}
This verifies that the $u^{GELDG} = u^{SLDG}$ for semi-discrete schemes.

  \end{proof}
\end{thm}
The fully discrete GEL RK DG and SL DG scheme are not equivalent for any $P^k$ approximation spaces.
The equivalence holds true for the special case of GEL RK DG $P^0$ and $P^1$ schemes, as specified in the Theorem below.
\begin{thm}
Under the same condition as Theorem \ref{lemma:equivalence between GELDG and SLDG}, the fully discrete GEL RK DG and SL DG scheme are equivalent for $P^k$ ($k\le1$) approximation spaces with any RK time discretization. Thus GEL RK DG schemes with $P^k$ ($k\le1$) approximation spaces are unconditionally stable.
  \begin{proof}
    We first consider the GEL DG method with forward-Euler time discretization by ~\eqref{the semi-discrete scheme for linear scaler equation}
    \begin{equation*}
    \begin{aligned}
      \int_{{I}_j}u^{GELDG}(x,t^{n+1})\Psi(x)dx&=\int_{\tilde{I}_j(t^n)}u^{n}(x)\Psi(x+(1+c)\Delta t)dx-c\Delta t\int_{\tilde{I}_j(t^n)} u^{n}(x) \Psi_x(x+(1+c)\Delta t)dx\\
      &\doteq RHS1^{GELDG}.
    \end{aligned}
%
    \end{equation*}
    We compute the right hand side
    \begin{equation*}
    \begin{aligned}
    &RHS^{SLDG}-RHS1^{GELDG}\\
    &=\int_{\tilde{I}_j(t^n)}u^{n}(x)(\Psi(x+\Delta t)-\Psi(x+(1+c)\Delta t))dx+\Delta t\int_{\tilde{I}_j(t^n)} c u^{n}(x) \Psi_x(x+(1+c)\Delta t)dx\\
    &=\int_{\tilde{I}_j(t^n)}u^{n}(x)(\Psi(x+\Delta t)-\Psi(x+(1+c)\Delta t)+c \Delta t \Psi_x(x+(1+c)\Delta t))dx
    \end{aligned}
    \end{equation*}
    So, if $\Psi(x)\in P^k, k\leq1$, then $RHS^{SLDG}-RHS1^{GELDG}=0$. That is, the GEL DG and SL DG scheme are equivalent for $P^k, k\leq1$ approximation space with forward-Euler time discretization.
    Such equivalence can be generalized to any SSP RK methods, which can be written as a convex combination of forward Euler method.
  \end{proof}
  \label{thm: fully_discrete_eqv}
\end{thm}
\begin{rem}
\label{rem: fully_discrete_eqv}
It can be shown that under the same assumption as Theorem~\ref{lemma:equivalence between GELDG and SLDG} with $P^2$ approximation space, the fully discrete GEL RK DG method, when the SSP RK2 and RK3 are applied for time discretization in Table \ref{table:ssprk}, is equivalent to the SL DG. Thus the scheme is unconditionally stable. However, the GEL RK DG, when coupled with forward Euler time discretization, is not equivalent to the SL DG and is not unconditionally stable.
%
\end{rem}
\subsection{Equivalence between semi-discrete GEL DG and EL DG scheme}\label{section:Equivalence between GELDG and ELDG}
\begin{thm}
\label{prop: equiv}
  For a linear transport problem with variable coefficient \eqref{scalar1d}, GEL DG scheme is equivalent to EL DG scheme in the semi-discrete form, assuming that they have the same space-time partition.
  \begin{proof}
  We first consider EL DG scheme for scaler equation \eqref{scalar1d}, which is formulated as
  \begin{equation}
\frac{d}{dt} \int_{\tilde{I}_j(t)}(u\psi)dx =-  \left(\hat{F} \psi\right) \left|_{\tilde{x}_{j+\frac12}(t) } \right. +   \left(\hat{F} \psi\right) \left|_{\tilde{x}_{j-\frac12}(t) } \right.   + \int_{\tilde{I}_j(t)}F\psi_xdx,
\label{mol_2ELDG}
\end{equation}
with $F(u) \doteq (a-\alpha) u$,
 \beq
\alpha(x, t)=-\nu_{j-\frac12} \frac{x-\tilde{x}_{j+\frac12}(t)}{\Delta x_j(t)} + \nu_{j+\frac12} \frac{x-\tilde{x}_{j-\frac12}(t)}{\Delta x_j(t)}\in P^1([\tilde{x}_{j-\frac12}(t), \tilde{x}_{j+\frac12}(t)]),
\label{eq: alpha_def1}
\eeq
and $\hat{F}$ is a monotone numerical flux.
We represent u in the form of \eqref{representation of basis}, and take $\psi$ as our basis function $ \tilde{\psi}_{j,m}(x,t)$ satisfy \eqref{eq: basis function_t} with $\Psi_{j,m}(x)$ being orthonormal basis of  the space $P^k(I_j)$, which follows the adjoint problem of EL DG:
\begin{equation}
\begin{cases}
\psi_t + \alpha(x,t) \psi_x =0 ,\ (x, t)\in\Omega_j,\\
\psi(t=t^{n+1}) = \Psi (x), \quad \forall \Psi(x) \in P^k(I_j).
\end{cases}
\label{eq: adjoint_of ELDG}
\end{equation}
Then we have $$\int_{\tilde{I}_j(t)} \tilde{\psi}_{j,i}(x,t) \tilde{\psi}_{j,m}(x,t)dx=\Delta x_j(t)\int_{I_j} \Psi_{j,i}(x)\Psi_{j,m}(x)dx=\delta_{i,m}\Delta x_j(t).$$  We can write the EL DG scheme \eqref{mol_2ELDG}
  into a matrix form
    \begin{equation}
    \frac{d}{dt} (U_j(t) \Delta x_j(t)) =-  \left(\hat{F} \bm{\tilde{\psi}_j}\right) \left|_{\tilde{x}_{j+\frac12}(t) } \right. +   \left(\hat{F} \bm{\tilde{\psi}_j}\right) \left|_{\tilde{x}_{j-\frac12}(t) } \right.   + \tilde{B}_j(t)U_j(t),
    \label{matrix form: semi-equivalence of ELDG}
    \end{equation}
    where $\tilde{B}_j(t)$ is a time-dependent matrix consisting of $[\tilde{B}_j(t)]_{i,m}=\int_{\tilde{I}_j(t)}(a(x,t)-\alpha)\tilde{\psi}_{j,m}(\tilde{\psi}_{j,i})_xdx, i=0,...k, m=0,...k$, and $\bm{\tilde{\psi}_j}=(\tilde{\psi}_{j,0},...,\tilde{\psi}_{j,k})^T$ is a vector of size $(k + 1)\times 1$.


   For GEL DG, we share the same space-time partition and we can rewrite the scheme \eqref{1D GELDG scheme for scalar} as
    \begin{equation}
    \frac{d}{dt} \int_{\tilde{I}_j(t)}(u\psi_{j,i}(x,t))dx =-  \left(\hat{F} \psi_{j,i}\right) \left|_{\tilde{x}_{j+\frac12}(t) } \right. +   \left(\hat{F} \psi_{j,i}\right) \left|_{\tilde{x}_{j-\frac12}(t) } \right.   + \int_{\tilde{I}_j(t)}(a(x,t)-\alpha)u(\psi_{j,i})_xdx,
    \label{semi-equivalence of GELDG}
    \end{equation}
    where $\psi_{j,i}$ is the test function as defined in \eqref{eq: test function_t}. Then, we rewrite the GEL DG scheme \eqref{semi-equivalence of GELDG} into the matrix form
    \begin{equation}
    \frac{d}{dt} (A_j(t) \Delta x_j(t) U_j(t)) =-  \left(\hat{F} \bm{{\psi}_j}\right) \left|_{\tilde{x}_{j+\frac12}(t) } \right. +   \left(\hat{F} \bm{{\psi}_j}\right) \left|_{\tilde{x}_{j-\frac12}(t) } \right.   + B_j(t)U_j(t),
    \label{matrix form: semi-equivalence of GELDG}
    \end{equation}
    where
    $[A_j(t)]_{i,m}=\int_{\tilde{I}_j(t)} \tilde{\psi}_{j,m}(x,t)\psi_{j,i}(x,t)dx /\Delta x_j(t)$, $[B_j(t)]_{i,m}=\int_{\tilde{I}_j(t)}(a(x,t)-\alpha_j)\tilde{\psi}_{j,m}(\psi_{j,i})_xdx$ and  $\bm{{\psi}_j}=({\psi}_{j,0},...,{\psi}_{j,k})^T$. Applying the product rule to the LHS of \eqref{matrix form: semi-equivalence of GELDG}, and with manipulations of the equation, we have
    \begin{equation}
    \frac{d}{dt} ( \Delta x_j(t) U_j(t)) =A_j^{-1}(t)\left(-  \left(\hat{F} \bm{{\psi}_j}\right) \left|_{j+\frac12} \right. +   \left(\hat{F} \bm{{\psi}_j}\right) \left|_{j+\frac12} \right.   + B_j(t)U_j(t)-\dot{A_j}(t)\Delta x_j(t) U_j(t)\right).
    \label{matrix form1: semi-equivalence of GELDG}
    \end{equation}
 As $\{\tilde{\psi}_{j,i}(x,t)\}_{i=0}^k$ and $\{\psi_{j,i}(x,t)\}_{i=0}^k$ are basis of $P^k(\tilde{I}_j(t))$, we can represent $\psi_{j,i}(x,t)$ as
 \[
 \psi_{j,i}(x,t)=\sum_{n=0}^k(\psi_{j,i}(x,t), \tilde{\psi}_{j,n}(x,t))_{\tilde{I}_j(t)}\tilde{\psi}_{j,n}(x,t)/\Delta x_j(t)=[A_j(t)]_{i}\bm{\tilde{\psi}_j},
 \]
 where $[A_j(t)]_{i}$ is the $i$th row of matrix $A_j(t)$.
    That is,
    \begin{equation}
      A_j^{-1}(t)\bm{{\psi}_j}=\bm{\tilde{\psi}_j},
    \label{transmision between basis function and test function}
    \end{equation}
    which leads to $$A_j^{-1}(t)(-  \left(\hat{F} \bm{{\psi}_j}\right) \left|_{j+\frac12} \right. +   \left(\hat{F} \bm{{\psi}_j}\right) \left|_{j+\frac12} \right.)=-  \left(\hat{F} \bm{\tilde{\psi}_j}\right) \left|_{\tilde{x}_{j+\frac12}(t) } \right. +   \left(\hat{F} \bm{\tilde{\psi}_j}\right) \left|_{\tilde{x}_{j-\frac12}(t) } \right..$$
    Next, we compute $\dot{A_j}(t)$, by $\Delta x_j(t)=\Delta x+(\nu_{j-\frac12}-\nu_{j+\frac12})(t^{n+1}-t)$ and $\partial_x(\alpha)=(\nu_{j+\frac12}-\nu_{j-\frac12})/\Delta x_j(t)$
    \begin{align*}
      &[\dot{A_j}(t)]_{i,m}=\frac{1}{\Delta x_j(t)}\frac{d}{dt}\int_{\tilde{I}_j(t)}\tilde{\psi}_{j,m}(x,t)\psi_{j,i}(x,t)dx
      +\frac{(\nu_{j-\frac12}-\nu_{j+\frac12})}{(\Delta x_j(t))^2}\int_{\tilde{I}_j(t)}\tilde{\psi}_{j,m}(x,t)\psi_{j,i}(x,t)dx\\
      &=\frac{1}{\Delta x_j(t)}\int_{\tilde{I}_j(t)}\left(\partial_t(\tilde{\psi}_{j,m})\psi_{j,i}+\tilde{\psi}_{j,m}\partial_t(\psi_{j,i})
      +\partial_x(\alpha\tilde{\psi}_{j,m}\psi_{j,i})\right)dx+\frac{(\nu_{j-\frac12}-\nu_{j+\frac12})}{(\Delta x_j(t))^2}\int_{\tilde{I}_j(t)}\tilde{\psi}_{j,m}(x,t)\psi_{j,i}(x,t)dx\\
      &=\frac{1}{\Delta x_j(t)}\int_{\tilde{I}_j(t)}\left(-\alpha \partial_x(\tilde{\psi}_{j,m})\psi_{j,i}-\tilde{\psi}_{j,m}\alpha_j\partial_x(\psi_{j,i})
      +\partial_x(\alpha)\tilde{\psi}_{j,m}\psi_{j,i}+\alpha\partial_x(\tilde{\psi}_{j,m})\psi_{j,i}+\alpha\tilde{\psi}_{j,m}\partial_x(\psi_{j,i})\right)dx\\
      &+\frac{(\nu_{j-\frac12}-\nu_{j+\frac12})}{(\Delta x_j(t))^2}\int_{\tilde{I}_j(t)}\tilde{\psi}_{j,m}(x,t)\psi_{j,i}(x,t)dx\\
      &=\frac{1}{\Delta x_j(t)}\left(\int_{\tilde{I}_j(t)} (\alpha-\alpha_j)\tilde{\psi}_{j,m}\partial_x(\psi_{j,i})dx+\frac{(\nu_{j+\frac12}-\nu_{j-\frac12})}{\Delta x_j(t)}\int_{\tilde{I}_j(t)}\tilde{\psi}_{j,m}\psi_{j,i}dx\right)\\
      &+\frac{(\nu_{j-\frac12}-\nu_{j+\frac12})}{(\Delta x_j(t))^2}\int_{\tilde{I}_j(t)}\tilde{\psi}_{j,m}(x,t)\psi_{j,i}(x,t)dx\\
      &=\frac{1}{\Delta x_j(t)}\int_{\tilde{I}_j(t)} (\alpha-\alpha_j)\tilde{\psi}_{j,m}(x,t)\partial_x(\psi_{j,i})(x,t)dx.
      \end{align*}
    So we have $[B_j(t)]_{i,m}-[\dot{A_j}(t)]_{i,m}\Delta x_j(t)=\int_{\tilde{I}_j(t)} (a(x,t)-\alpha)\tilde{\psi}_{j,m}(x,t)\partial_x(\psi_{j,i})(x,t)dx$. Then we can easily get $A_j^{-1}(t)(B_j(t)-\dot{A_j}(t)\Delta x_j(t))=\tilde{B}_j(t)$ by \eqref{transmision between basis function and test function}, which shows the equivalence of \eqref{matrix form: semi-equivalence of ELDG} and \eqref{matrix form1: semi-equivalence of GELDG}.
  \end{proof}
\end{thm}

\begin{rem}
 For general nonlinear problems, the equivalence of semi-discrete GEL DG and EL DG methods can be established in a similar way, given the same space-time partition.
\end{rem}

\begin{rem} (Fully discrete case)
The fully discrete EL RK DG scheme is known as equivalent to the fully discrete ALE DG scheme combined with one extra step of solution projection (which does not affect stability). Thus the stability result of ALE DG method for linear conservation laws \cite{Zhou2019arbitrary} can be directly applied to assess the stability property of fully discrete EL RK DG scheme. The stability of fully discrete GEL RK DG method is still theoretically open, and is being investigated numerically in the numerical section.
\end{rem}

\section{DGCL, MPP properties and numerical limiters}
\label{section:GCL and MPP}

We consider the property of DGCL, which requires that the numerical scheme reproduces exactly a constant solution under the geometric parameters of numerical schemes. As shown in \cite{CharbelFarhat2001GCL}, satisfaction of DGCL is a necessary and sufficient condition for a numerical scheme to preserve the nonlinear stability of its fixed grid counterpart. In the context of GEL DG scheme, the geometric parameters refer to the parameters involved in the space-time partition on which the PDE is evolved. It was established in \cite{klingenberg2017arbitrary} that ALE-DG scheme (the evolution step in the EL DG method) satisfies the DGCL for 1D problems, and for high D problems with the time integrator which holds the accuracy not less than the value of the spatial dimension \cite{Fu2019arbitrary}. Below we shown the conditions under which the DGCL holds for GEL DG method when coupled with forward Euler time discretization, see Table~\ref{The GCL related parameters settings}. As a direct consequence, we find the DGCL no longer holds for general GEL RK DG schemes. Proposition~\ref{prop: DGCL} is numerically verified in Table~\ref{Check GCL numerically}.
\begin{table}[!ht]\footnotesize
\caption{The related parameters settings about RHS and DGCL of GEL DG method for $u_t+au_x=0.$
  }
\centering
\begin{tabular}{| c | c | c|   }

 \hline
 $\Psi$ &  RHS & condition for DGCL
  \\
  \hline
      $1$ & $\Delta x$ &  -- \\
  \hline
    $x$ & $\frac{\Delta t}{\Delta x}(\nu_{j-\frac12}-\nu_{j+\frac12})+\frac{\Delta t^2}{\Delta x}(\nu_{j-\frac12}-\nu_{j+\frac12})(\frac{\nu_{j-\frac12}+\nu_{j+\frac12}}{2}-\alpha_j)$        &     $\nu_{j+\frac12}=\nu_{j-\frac12}$         \\
    \hline
     $x^2-\frac{1}{12}$ &  $-\frac{\Delta t^2}{2\Delta x}[(\nu_{j+\frac12}-\alpha_j)^2+(\nu_{j-\frac12}-\alpha_j)^2 ]+\frac{2\Delta t^3}{3\Delta x^2}[(\nu_{j+\frac12}-\alpha_j)^3-(\nu_{j-\frac12}-\alpha_j)^3] $ &
     $\nu_{j+\frac12}=\nu_{j-\frac12}=\alpha_j$      \\
\hline
  \end{tabular}
\label{The GCL related parameters settings}
\end{table}
\begin{table}[!ht]\footnotesize
\caption {DGCL. The error $\|u_h(x,t)-u(x,t)\|_{L^\infty}$ of the GEL DG scheme solving $u_t+u_x=0$ with initial condition $u(x,0) = 1$, $T=1$, $CFL=0.1$. We test schemes with different choices of $\nu_{j\pm\frac12}$ as parameters for space-time partition and $\alpha_j$ as parameters for adjoint problems.}
\centering
\begin{tabular}{| c| c| c| c| c| }
\hline
\diagbox{\ \ $\nu_{j\pm\frac{1}{2}}$\ \ }{\ $\alpha_j$\ \ \ } & & $1$  &  $1+0.5$ & $1+\Delta x\sin(x_j)$  \\
\hline
       &  $P^0$ &   $10^{-16}$  &   $10^{-16}$ &    $10^{-16}$     \\
    1  &  $P^1$ &   $10^{-15}$  &   $10^{-15}$ &    $10^{-15}$     \\
       &  $P^2$ &   $10^{-15}$  &   $10^{-02}$  &   $10^{-02}$      \\
\hline
       &  $P^0$ &   $10^{-15}$  &   $10^{-15}$ &    $10^{-15}$     \\
 1+0.5 &  $P^1$ &   $10^{-15}$  &   $10^{-15}$ &    $10^{-15}$     \\
       &  $P^2$ &   $10^{-02}$   &   $10^{-14}$ &   $10^{-02}$      \\
\hline
       &  $P^0$ &   $10^{-15}$  &   $10^{-15}$ &    $10^{-15}$     \\
 $1+\Delta x\sin(x_{j\pm\frac{1}{2}})$  &  $P^1$ &   $10^{-03}$  &   $10^{-03}$ &    $10^{-05}$     \\
       &  $P^2$ &   $10^{-03}$  &   $10^{-02}$  &   $10^{-04}$      \\
\hline
\end{tabular}
\label{Check GCL numerically}
\end{table}
\begin{prop}
\label{prop: DGCL}
Under the conditions specified in Table~\ref{The GCL related parameters settings} that the GEL DG method coupled with forward Euler time discretization satisfies the DGCL.
\begin{proof}
    For the GEL DG scheme with Forward-Euler time discretization, we can get
\begin{equation}\label{eqn: forward euler time discretization}
\begin{aligned}
  &\int_{\tilde{I}_j}(u_j^{n+1}\psi(x,t^{n+1}))dx =\int_{I^*_j}(u_j^n\psi(x,t^n))dx +\Delta t(-  \left(\hat{F}\psi\right) \left|_{\tilde{x}_{j+\frac12}(t^n) } \right. +   \left(\hat{F}\psi\right) \left|_{\tilde{x}_{j-\frac12}(t^n) } \right.)\\
  &   + \Delta t \int_{I^*_j} (a-\alpha_j)u \psi_xdx\doteq RHS.
\end{aligned}
\end{equation}
We let $u_j^n=1$,
\begin{align*}
  RHS&=\int_{I^*_j} \Psi(x+\alpha_j \Delta t)dx+\Delta t(-(a-\nu)\Psi(x+\alpha_j \Delta t)|_{\tilde{x}_{j+\frac12}(t^n) }+(a-\nu)\Psi(x+\alpha_j \Delta t)|_{\tilde{x}_{j-\frac12}(t^n) }\\
  &+(a-\alpha_j)\Psi(x+\alpha_j \Delta t)|_{\tilde{x}_{j+\frac12}(t^n) }-(a-\alpha_j)\Psi(x+\alpha_j \Delta t)|_{\tilde{x}_{j-\frac12}(t^n)})\\
  &=\int_{I^*_j} \Psi(x+\alpha_j \Delta t)dx+\Delta t(-(\alpha_j-\nu)\Psi(x+\alpha_j \Delta t)|_{\tilde{x}_{j+\frac12}(t^n) }+(\alpha_j-\nu)\Psi(x+\alpha_j \Delta t)|_{\tilde{x}_{j-\frac12}(t^n) }).
\end{align*}
With different choices of $\Psi$, we show conditions in Table \ref{The GCL related parameters settings}, under which the DGCL is satisfied.
  \end{proof}
\end{prop}
From the above proposition, the GEL DG method with $P^0$ polynomial space and forward Euler satisfies DGCL. However, for the method with high order polynomial spaces, DGCL fails for general GEL RK DG schemes. On the other hand, MPP is the principle that numerical solutions will be bounded by the maximum and minimum of initial condition. One can show that if a scheme satisfies the MPP, it automatically satisfies the DGCL.
 Below, we are going to first show that first order GEL DG scheme satisfies the MPP property.
Then we apply two MPP limiters to high order GEL RK DG schemes to get the MPP, hence the DGCL property.
\begin{prop} (MPP property of the first order GEL DG scheme.)
Let $u_m=\min{u_0(x)}, u_M=\max{u_0(x)}$, then the first order GEL DG solution $\bar{u}_j^{n}\in [u_m,u_M]$, $\forall j, n$ under the condition $\Delta t\leq \frac{\Delta x}{\alpha_1}$ where $\alpha_1=\max_j|a-\nu|_{j+ \frac12}$.
\begin{proof}
  The first order GEL DG scheme reads
  \begin{equation*}
     \int_{\tilde{I}_j}\bar{u}_j^{n+1}dx =\int_{I^*_j}\bar{\tilde{u}}_j^{n}dx -\Delta t\left( \hat{h}_{j+ \frac12}-\hat{h}_{j- \frac12}  \right)
  \end{equation*}
  where $\hat{h}$ is the first order monotone flux, e.g. the Lax-Friedrich flux,
  \[
  \hat{h}_{j+ \frac12}=\hat{F}|_{\tilde{x}_{j+\frac12}(t^n) } =\frac{(a-\nu)_{j+ \frac12}(\bar{\tilde{u}}_{j+1}^{n}+\bar{\tilde{u}}_j^{n})}{2}-\frac{\alpha_1(\bar{\tilde{u}}_{j+1}^{n}-\bar{\tilde{u}}_j^{n})}{2},
 \]
 and $\tilde{u}_j^{n}$ is the reaveraging  of solution ${u}_j^{n}$ on $I_j^*$.
$\bar{\tilde{u}}_j^{n}\in [u_m,u_M]$, since ${u}_j^{n}\in [u_m,u_M]$.
Let $\lambda=\frac{\Delta t}{\Delta x}$, we get
\begin{align*}
  \bar{u}_j^{n+1} &=\frac{\Delta x_j^*}{\Delta x}\bar{\tilde{u}}_j^{n}- \lambda [\hat{h}_{j+ \frac12}-\hat{h}_{j- \frac12}]\\
   &=(1-\lambda \alpha_1)\bar{\tilde{u}}_j^{n}+\lambda\frac{\alpha_1-(a-\nu)_{j+ \frac12}}{2}\bar{\tilde{u}}_{j+1}^{n}+\lambda\frac{(a-\nu)_{j- \frac12}+\alpha_1}{2}\bar{\tilde{u}}_{j-1}^{n}.
\end{align*}
Under the condition of $\Delta t\leq \frac{\Delta x}{\alpha_1}$ with $\alpha_1=\max_j|a-\nu|_{j+ \frac12}$, the coefficients $1-\lambda \alpha_1, \alpha_1\pm(a-\nu)_{j- \frac12}$ are all positive, hence $\bar{u}_j^{n+1}$ is a convex combination of $\bar{\tilde{u}}_j^{n},\bar{\tilde{u}}_{j+1}^{n},\bar{\tilde{u}}_{j-1}^{n}$. Therefore,
$\bar{u}_j^{n+1} \in[u_m,u_M]$.
\end{proof}
\label{prop: mpp}
\end{prop}

Next we propose to apply MPP limiters to high order GEL RK DG schemes to obtain the MPP property, leading to DGCL in a general setting. There are two MPP limiters that we will make use of: one is the rescaling limiter for DG polynomials \cite{Zhang2010MPPlimiter, ZSpos}, and the other is the parametrized flux limiter in a sequence of work \cite{xiong2013parametrized, xiong2014high, Xiong2015MPPlimiter}. While the rescaling limiter works very well for RKDG schemes, it often leads to extra time step restrictions for MPP property and reduction in temporal order of convergence \cite{Zhang2010MPPlimiter}. For the GEL RK DG scheme, we propose to apply the rescaling limiter \cite{Zhang2010MPPlimiter, ZSpos} only for solutions at $t^n$ and in the solution remapping step, and apply the flux limiter only in the final stage of RK methods \cite{xiong2013parametrized}. Neither the rescaling limiter nor the flux limiter are applied in internal stages of RK methods to avoid the temporal order reduction of RK methods.

Below, we
demonstrate the procedure of the proposed MPP flux limiter in the context of a third order SSP RK time discretization.

Step 1:  We first use the polynomial rescaling MPP limiter \cite{Zhang2010MPPlimiter} to preserve the MPP property for the DG polynomial $u(x, t^n)$ on $I_j$, then perform the $L^2$ projection of $u(x, t^n)$ from background cells onto upstream cells ${I}^*_j$. We assume $\tilde{u}_j^{n}$ as the reprojected polynomial on $I_j^*$, then $\bar{\tilde{u}}_j^{n}\in [u_m,u_M]$.

Step 2: We propose to apply the parametrized MPP flux limiters \cite{Xiong2015MPPlimiter} in the context of a moving mesh to guarantee $\bar{u}_j^{n+1}\in [u_m,u_M]$.
Taking $\psi = 1$ in \eqref{1D GELDG scheme for scalar}, we have
\begin{equation}\label{semi-discrete SSPRK3 MPP}
     \frac{d}{dt}\int_{\tilde{I}_j(t)}u_j(x,t)dx = \hat{H}_{j+ \frac12}-\hat{H
}_{j- \frac12},
\end{equation}
where $ \hat{H}$ is the high order monotone flux. With the third order SSP RK method, the
update of cell averages in equation \eqref{semi-discrete SSPRK3 MPP} can be written as
\begin{equation*}
\bar{u}_j^{n+1} =\frac{\Delta x_j^*}{\Delta x}\bar{\tilde{u}}_j^{n}- \lambda [\hat{H}^{rk}_{j+ \frac12}-\hat{H}^{rk}_{j- \frac12}],
\end{equation*}
where $\hat{H}^{rk}_{j+ \frac12}=\frac{1}{6}\hat{H}^{n}_{j+ \frac12}+\frac{1}{6}\hat{H}^{(1)}_{j+ \frac12}+\frac{2}{3}\hat{H}^{(2)}_{j+ \frac12}$ with $\hat{H}^{(i)},i=0,1,2$ being the numerical flux obtained by $u_h(x,t^{(i)})$ at each RK stage. The MPP flux limiter is proposed to replace the numerical flux $\hat{H}^{rk}_{j+ \frac12}$
by a modified one $\tilde{H}^{rk}_{j+ \frac12}=\theta_{j+ \frac12}(\hat{H}^{rk}_{j+ \frac12}-\hat{h}_{j+ \frac12})+\hat{h}_{j+ \frac12}$, where $\hat{h}_{j+ \frac12}$ is the numerical flux from a first order scheme as in Proposition~\ref{prop: mpp} and $\theta_{j+ \frac12} \in [0, 1]$, $\forall j$.  The parameter $\theta_{j+ \frac12}$ is set to be as close to $1$ as possible, and to ensure $\bar{u}_j^{n+1}\in [u_m,u_M]$, for which sufficient inequalities to have are, $\forall j $,
\begin{align}
  \lambda \theta_{j- \frac12}(\hat{H}^{rk}_{j- \frac12}-\hat{h}_{j- \frac12})-\lambda \theta_{j+ \frac12}(\hat{H}^{rk}_{j+ \frac12}-\hat{h}_{j+ \frac12})-\Gamma_j^M&\leq0, \label{xulimiterM}\\
  \lambda \theta_{j- \frac12}(\hat{H}^{rk}_{j- \frac12}-\hat{h}_{j- \frac12})-\lambda \theta_{j+ \frac12}(\hat{H}^{rk}_{j+ \frac12}-\hat{h}_{j+ \frac12})-\Gamma_j^m&\geq0, \label{xulimiterm}
\end{align}
with $$\Gamma_j^M= u_M-\frac{\Delta x_j^*}{\Delta x}\bar{u}_j^n+\lambda (\hat{h}_{j+ \frac12}-\hat{h}_{j- \frac12}),\ \Gamma_j^m= u_m-\frac{\Delta x_j^*}{\Delta x}\bar{u}_j^n+\lambda (\hat{h}_{j+ \frac12}-\hat{h}_{j- \frac12}).$$ Here we adjust the terms $\Gamma_j^M$ and $\Gamma_j^m$ from \cite{Xiong2015MPPlimiter} in the context of moving meshes, and obtain the parameter $\theta_{j+ \frac12}$ satisfying ~\eqref{xulimiterM} and ~\eqref{xulimiterm} for all j, thus guarantee $u_h^{n+1}\in [u_m,u_M]$. Note that such $\theta$ always exists, since $\theta_{j+\frac12}=0$ is a solution to  ~\eqref{xulimiterM} and ~\eqref{xulimiterm} for all $j$. Then we go back to Step 1 to apply the polynomial rescaling MPP limiter again to ensure DG polynomials $u_h(x,t^{n+1})$ satisfy the MPP property.


We call the GEL RK DG method with the above described limiter `GELDGMPPlimiter'. We also apply the polynomial rescaling MPP limiter \cite{Zhang2010MPPlimiter} at each RK internal stages, for which the scheme is termed `zhangMPPlimiter'. We will compare numerical performance, in terms of error and time stepping size for numerical stability, of these two limiters in the next section.
\begin{rem}
 For a linear variable coefficient equation, MPP property is lost, but the PP property stays valid. A PP limiter can be applied in a similar fashion to preserve the PP property.
 \end{rem}

\section{Numerical results}
\label{section:numerical}

In this section, we perform numerical experiments for linear transport problems, where we set the time stepping size as $\Delta t= \frac{CFL}{a} \Delta x$ for 1D and $\Delta t=\frac{CFL}{\frac{a}{\Delta x}+\frac{b}{\Delta y}}$, where a and b are maximum transport speed in $x$- and $y$-directions respectively. We mainly study the following aspects: the spatial order of convergence by using small enough time stepping size, the temporal order of convergence by varying CFL, numerical stability under a large time stepping size. When applicable, we also present the EL DG solutions \cite{cai2020eldg} for comparison.


\subsection{1D linear transport problems}
\begin{exa}
(1D linear transport equation with constant coefficient.)
We consider a simple 1D transport equation
\begin{equation}\label{eqn:transport with constant coefficient}
u_t + u_x =0, \ x\in[0,2\pi],
\end{equation}
with the smooth initial data $u(x,0)=\sin(x)$ and exact solution $u(x,t)=\sin(x-t)$.
For the constant coefficient problem, the proposed GEL DG method, if using the exact velocity field for space-time partition and the adjoint problem, is the same as EL DG and SL DG.
Here we perturb the velocity at cell boundaries, i.e. $\nu_{j+\frac12}$ and the velocity in the modified adjoint problem i.e. $\alpha_j$ in \eqref{eq: adjoint_new} to get GELDG1, GELDG2 and GELDG3 schemes respectively. Parameters of these GEL DG methods are given in Table \ref{The related parameters settings}.
\begin{table}[!ht]\footnotesize
\caption{The numerical parameters of ELDG, GELDG1, GELDG2 and GELDG3 method for ~\eqref{eqn:transport with constant coefficient}.
  }
\centering
\begin{tabular}{| c | c | c| c| c| c| }

 \hline
 &  \multicolumn{1}{c|}{ ELDG} & \multicolumn{1}{c|}{ GELDG1}& \multicolumn{1}{c|}{ GELDG2}& \multicolumn{1}{c|}{ GELDG3}
  \\ \hline
    mesh $\nu_{j+\frac12}$ &  $1+\sin(x_{j+\frac12})\Delta x$   &  $1+\sin(x_{j+\frac12})\Delta x $  &   1 &  $1+\sin(x_{j+\frac12})\Delta x $\\
  \hline
    adjoint $\alpha_j$     &          --         &     1       &     $1+\sin(x_j)\Delta x$    &     $1+\sin(x_j)\Delta x $   \\
\hline
  \end{tabular}
\label{The related parameters settings}
\end{table}
Table \ref{linear1d_spatial} reports the spatial accuracies of these methods for this example with the same time stepping size. The proposed GEL DG methods are found to perform comparably as the ELDG method. We vary time stepping size, with fixed well-resolved spatial meshes, and plot error vs. $CFL$ in Figure \ref{linear_stability} for these schemes with $P^1$-SSP RK2 solutions (left) and $P^2$-SSP RK3 solutions  (right) at a long integration time $T=100$. GELDG2 methods are found to be unconditionally stable with the space-time partition exactly following the characteristics, which is consistent with Theorem~\ref{thm: fully_discrete_eqv} and Remark~\ref{rem: fully_discrete_eqv}. GELDG3 behaves closer to ELDG and performs better than GELDG1, as the adjoint problem and the space-time partition are more closely related. It indicates that, designing the GEL DG scheme associating the adjoint problem with the space-time partition is advantageous for better performance of the scheme. In particular, the adjoint problem is determined by the space-time partition in the EL DG algorithm for which best numerical stability is observed and theoretically investigated in \cite{cai2020eldg}.

Further, we apply higher order RK methods for time discretization, as we are interested in using relatively large time stepping size. We show error vs. $CFL$ in Figure \ref{linear_stability1} for GELDG1 and GELDG3 schemes $P^1$ with SSP RK2, RK4 (left) and $P^2$ SSP RK3, RK4 (right) at a long integration time $T=100$. We can concluded that GELDG3 is better than GELDG1 in terms of stability and the higher order RK help with reducing the error magnitute when large time stepping size is used. We note that in both Figure \ref{linear_stability} and \ref{linear_stability1}, the CFL allowed with stability is much larger than that of the RK DG method which is $\frac{1}{2k+1}$.

\begin{table}[!ht]\footnotesize
\caption{1D linear transport equation with constant coefficient. $u_t+u_x=0$ with initial condition $u(x,0) = \sin(x)$. $T=\pi$.
We use $CFL=0.3$ and $CFL=0.18$ for all $P^1$ and $P^2$ schemes, respectively.}
\centering
\begin{tabular}{| c | cc  | cc| cc| cc| cc| }

\hline
Mesh  &{$L^1$ error} & Order  &  {$L^1$ error} & Order &  {$L^1$ error} & Order  &  {$L^1$ error} & Order \\
 \hline
  &  \multicolumn{2}{c|}{$P^1$ ELDG} & \multicolumn{2}{c|}{$P^1$ GELDG1}& \multicolumn{2}{c|}{$P^1$ GELDG2}& \multicolumn{2}{c|}{$P^1$ GELDG3}
  \\ \hline
    40 &     6.08E-04 &    --     &    6.08E-04 &    --    &    6.37E-04 &    --   &    6.08E-04 &    --   \\
    80 &     1.55E-04 &     1.97  &    1.55E-04 &    1.97  &    1.59E-04 &    2.00 &    1.55E-04 &    1.97   \\
   160 &     3.84E-05 &     2.02  &    3.84E-05 &    2.02  &    3.90E-05 &   2.03  &    3.84E-05 &    2.02   \\
   320 &     9.77E-06 &     1.98  &    9.77E-06 &    1.98  &    9.83E-06 &    1.99  &   9.77E-06 &    1.98   \\
\hline
  &  \multicolumn{2}{c|}{$P^2$ ELDG} & \multicolumn{2}{c|}{$P^2$ GELDG1}& \multicolumn{2}{c|}{$P^2$ GELDG2}& \multicolumn{2}{c|}{$P^2$ GELDG3}
  \\ \hline
    40 &     7.69E-06 & --      &    2.60E-05 &    --    &    7.25E-06  &    --   &    7.71E-06  &    --    \\
    80 &     9.45E-07 &     3.03  &    1.91E-06 &    3.77  &    9.23E-07 &   2.97 &    9.45E-07  &    3.03   \\
   160 &    1.18E-07 &     3.00  &    1.67E-07 &   3.51   &    1.17E-07 &    2.98  &    1.18E-07  &   3.00   \\
   320 &    1.41E-08 &     3.07   &    1.63E-08 &   3.36  &    1.40E-08 &    3.06  &    1.41E-08  &   3.07  \\
\hline
\end{tabular}
\label{linear1d_spatial}
\end{table}
\begin{figure}[!ht]
\centering
\includegraphics[height=50mm]{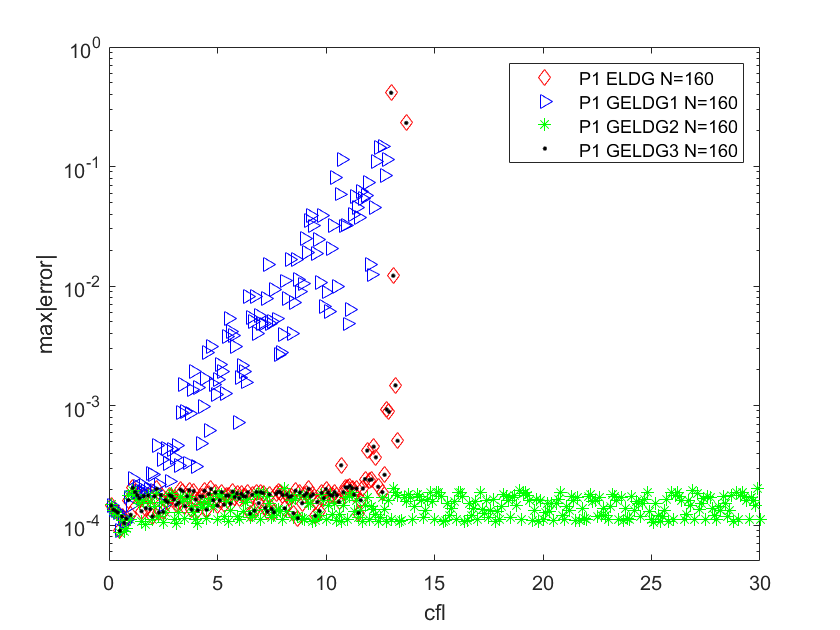}
\includegraphics[height=50mm]{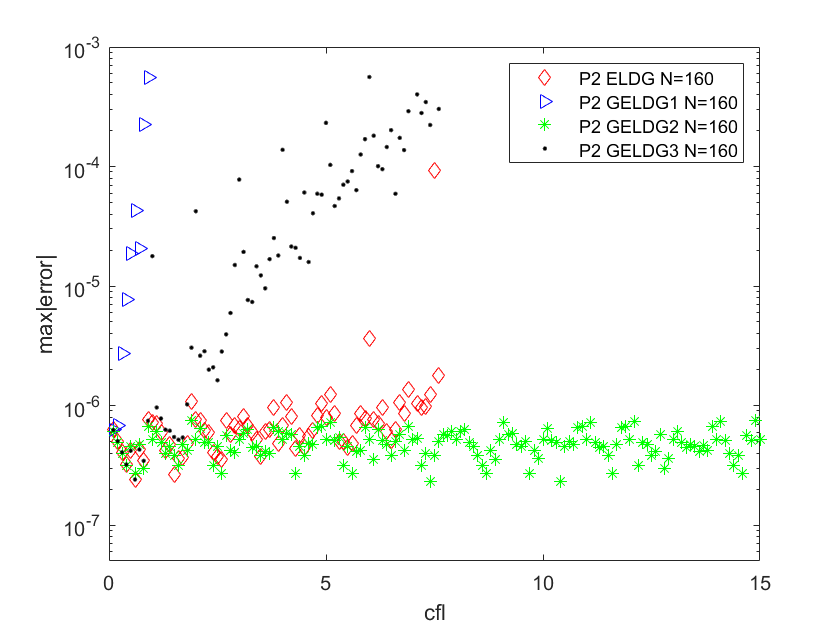}
\caption{The $L^\infty$ error versus $CFL$ of ELDG methods, GELDG1, GELDG2 and GELDG3 methods $P^1$ (left) with SSP RK2 and $P^2$ (right) with SSP RK3 time discretization for ~\eqref{eqn:transport with constant coefficient} with initial condition $u(x,0) = \sin(x)$. A long time simulation is performed with $T=100$ and mesh size $N=160$.  
}
\label{linear_stability}
\end{figure}

\begin{figure}[!ht]
\centering
\includegraphics[height=50mm]{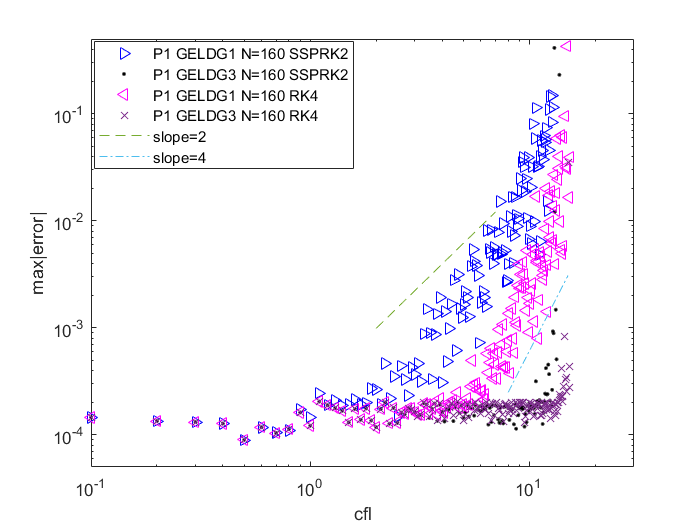}
\includegraphics[height=50mm]{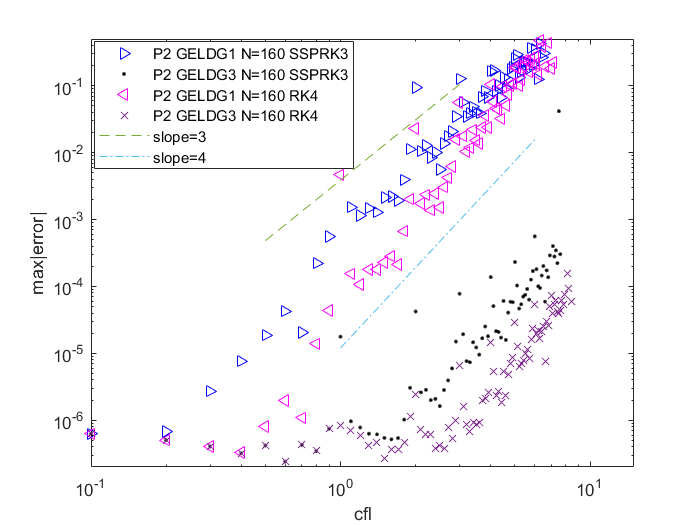}
\caption{The $L^\infty$ error versus $CFL$ of GELDG1 and GELDG3 methods $P^1$ (left) with SSP RK2, RK4 and $P^2$ (right) schemes SSP RK3, RK4 time discretization for ~\eqref{eqn:transport with constant coefficient} with initial condition $u(x,0) = \sin(x)$. A long time simulation is performed with $T=100$ and mesh size $N=160$.  
}
\label{linear_stability1}
\end{figure}

Finally, we verify the DGCL property of the scheme when the proposed MPP limiters are applied.
Table \ref{Check GCL numerically after MPP} shows that we regain the DGCL property by applying `GELDGMPPlimiter' to GEL RK DG schemes, whose DGCL property are assessed in Table \ref{Check GCL numerically}.
Next, we present in Figure \ref{MPP limiter} the $L^\infty$ error versus $CFL$ of GEL DG methods with `GELDG MPP limiter' and `zhang MPP limiter' $P^1$ and $P^2$ schemes. We apply the RK4 for time integration. From Figure \ref{MPP limiter}, we can observe better stability and accuracy for GEL DG method with `GELDGMPPlimiter' which is applied only in the final RK stage, compared with one with `zhangMPPlimiter' which is applied in every RK intermediate stages, especially for schemes with $P^2$ polynomial space.
We also test the MPP property for a step function initial condition:
\begin{equation}\label{eqn: initail step function}
  u(x,0)=\left\{
  \begin{aligned}
    &1, \ &2<x<7,\\
    &0, \ &otherwise.
  \end{aligned}
  \right.
\end{equation}
The computational domain is [0,90]. The solution of GELDG2 method for $P^2$ with $CFL=1$, $N=160$ and with RK4 time discretization without (left) and with (right) `GELDGMPPlimiter' are shown in Figure \ref{step function solution MPP}. We can observe the MPP property for GEL DG with `GELDGMPPlimiter'.

\begin{table}[!ht]\footnotesize
\caption{Compute $u_h(x,t)-u(x,t)$ to check if all $P^0$, $P^1$ and $P^2$ GEL DG schemes with `GELDGMPPlimiter' satisfy DGCL for ~\eqref{eqn:transport with constant coefficient} with initial condition $u(x,0) = 1$, $T=1$, $CFL=0.1$ respectively.}
\centering
\begin{tabular}{| c| c| c| c| c| }
\hline
\diagbox{\ \ $\nu_{j\pm\frac{1}{2}}$\ \ }{\ $\alpha_j$\ \ \ } & & $1$  &  $1+0.5$ & $1+\Delta x\sin(x_j)$  \\
\hline
       &  $P^0$ &   $10^{-16}$  &   $10^{-16}$ &   $10^{-16}$     \\
    1  &  $P^1$ &   $10^{-16}$  &   $10^{-16}$ &   $10^{-16}$     \\
       &  $P^2$ &   $10^{-16}$  &   $10^{-15}$ &   $10^{-16}$      \\
\hline
       &  $P^0$ &   $10^{-16}$  &   $10^{-16}$ &    $10^{-16}$    \\
 1+0.5 &  $P^1$ &   $10^{-16}$  &   $10^{-16}$ &    $10^{-16}$     \\
       &  $P^2$ &   $10^{-16}$  &   $10^{-15}$ &    $10^{-16}$      \\
\hline
       &  $P^0$ &   $10^{-16}$  &   $10^{-16}$ &    $10^{-16}$     \\
 $1+\Delta x\sin(x_{j\pm\frac{1}{2}})$  &  $P^1$ &  $10^{-16}$  &  $10^{-16}$ &    $10^{-16}$     \\
       &  $P^2$ &   $10^{-16}$  &   $10^{-15}$  &    $10^{-16}$      \\
\hline
\end{tabular}
\label{Check GCL numerically after MPP}
\end{table}

\begin{figure}[!ht]
\centering
\includegraphics[height=50mm]{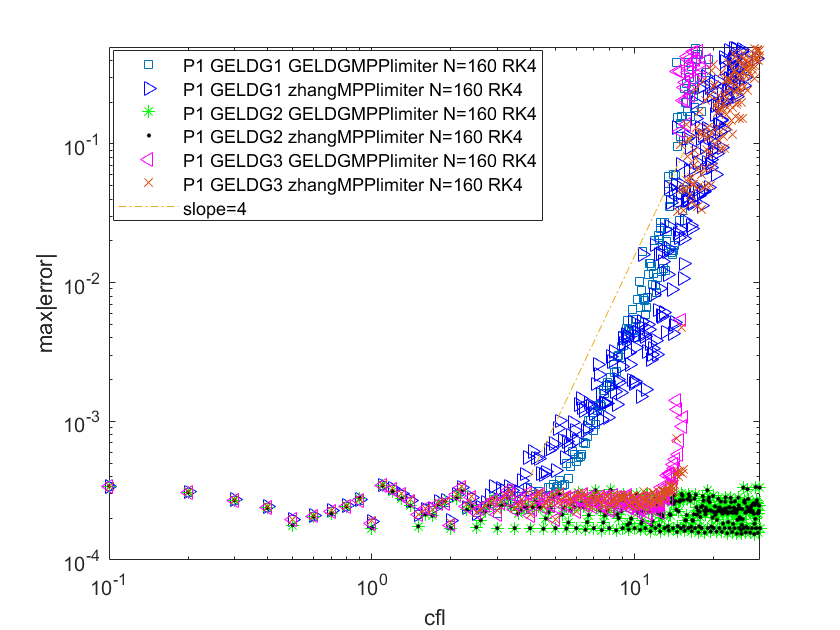}
\includegraphics[height=50mm]{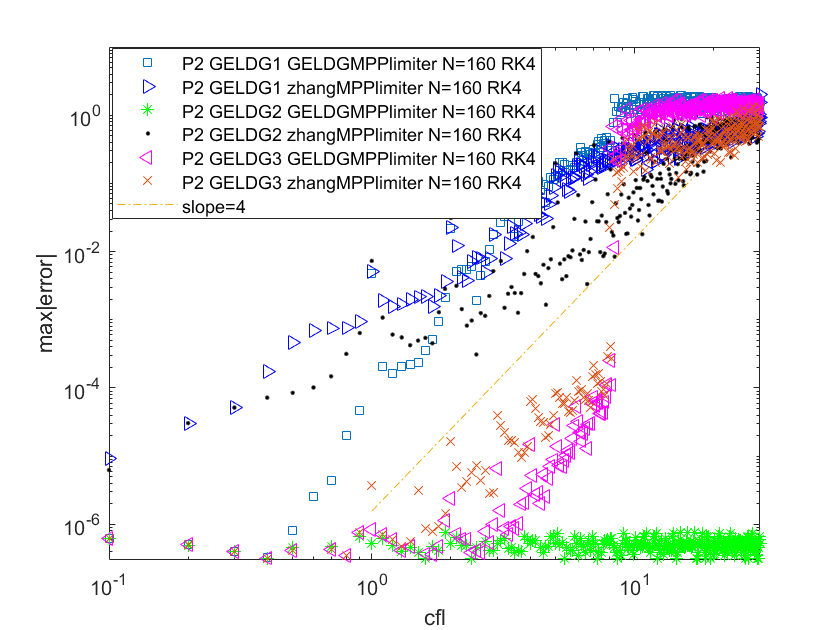}
\caption{The $L^\infty$ error versus $CFL$ of GELDG1, GELDG2 and GELDG3 methods with `GELDG MPP limiter' and `zhang MPP limiter' for $P^1$ (left) and $P^2$ (right) with RK4 time discretization for ~\eqref{eqn:transport with constant coefficient} with initial condition $u(x,0) = \sin(x)$. A long time simulation is performed with $T=100$ and mesh size $N=160$.  
}
\label{MPP limiter}
\end{figure}
\begin{figure}[!ht]
\centering
\includegraphics[height=50mm]{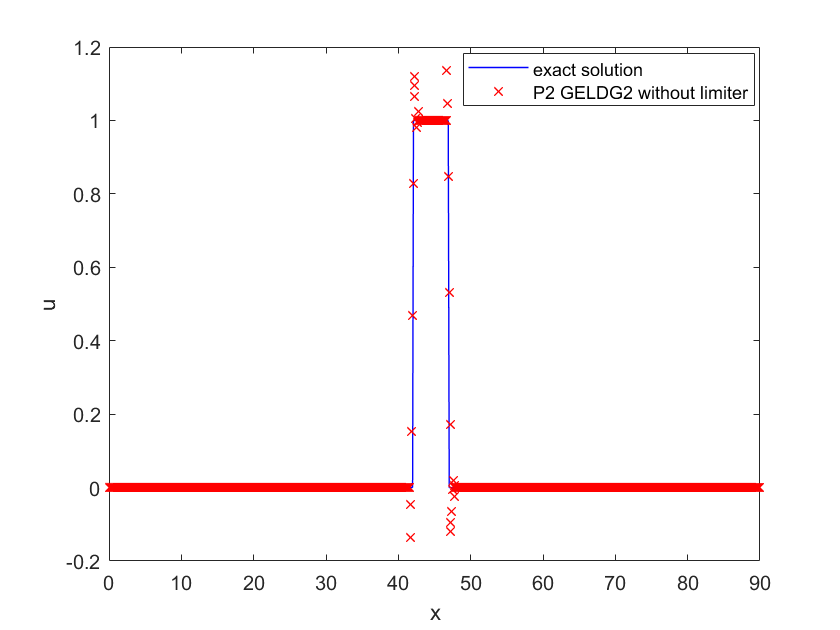}
\includegraphics[height=50mm]{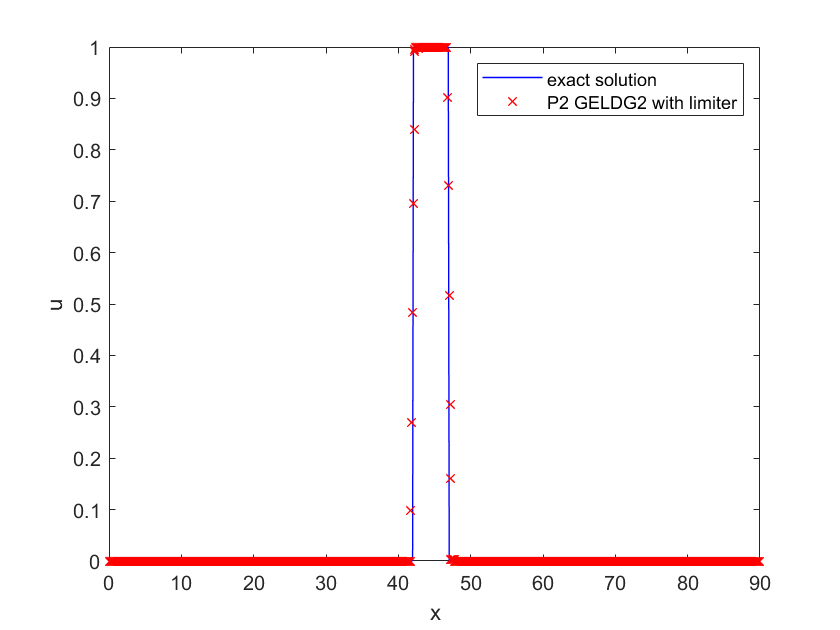}
\caption{The solution u of GELDG2 method $P^2$ without (left) and with (right) `GELDGMPPlimiter' with RK4 time discretization and $CFL=1$ for \eqref{eqn:transport with constant coefficient}
with initial condition \eqref{eqn: initail step function}. A long time simulation is performed with $T=40$ and mesh size $N=160$.  
}
\label{step function solution MPP}
\end{figure}

\end{exa}

\begin{exa}
(1D transport equation with variable coefficients.)
Consider
\begin{equation}\label{eqn: 1D transport equation with variable coefficients}
u_t + (\sin(x)u)_x =0, \ x\in[0,2\pi]
\end{equation}
with initial condition $u(x,0)=1$ and the periodic boundary condition.
The exact solution is given by
\begin{equation}
u(x,t) =
\frac{  \sin( 2\tan^{-1}( e^{-t}\tan(\frac{x}{2}) ) )   }{  \sin(x) }.
\end{equation}
The related parameters settings are given in Table \ref{The related parameters settings for example2}. The expected spatial convergence of ELDG and GELDG are shown in Table \ref{1dsin_example}. In Figure \ref{sin1d_time}, we plot the $L^\infty$ error versus $CFL$ of ELDG and GELDG scheme with $P^1$ (left) and $P^2$ (right) polynomial spaces. The following observations are made:
(1) both methods perform similarly around and before $CFL=1$, which is well above the stability constraint of the RK DG method $1/(2k+1)$; (2) after $CFL=1$ and before stability constraint of the method, the temporal convergence order is observed to be consistent with the order of RK discretization; (3) EL RK DG has better performance than GEL RK DG with the same mesh.

In addition, we test the proposed inflow boundary condition for the following problem
\begin{equation}
\begin{cases}
u_t + ( \sin(x)u)_x = 0, x\in [\frac{\pi}{2},\frac{5\pi}{2}]\\
u(x,0) = 1,\\
u(\frac{\pi}{2},t) =   \sin( 2\tan^{-1}( e^{-t}\tan(\frac{\pi}{4}) ) ).
\end{cases}
\label{eqn:1d}
\end{equation}
We use the GEL DG method with PP limiter to solve this problem, the $L^1$ and $L^\infty$ errors for $P^1$ and $P^2$ are shown in Table \ref{1dsin_example inflow boundary}. The optimal convergence rate are observed. Besides, we show the $L^\infty$ error versus $CFL$ of  GEL DG method
with PP limiter in Figure \ref{sin1d_time inflow boundary}. Expected temporal convergence is observed.

\begin{table}[!ht]\footnotesize
\caption{The related parameters settings of ELDG, GELDG method for $u_t+(\sin(x)u)_x=0.$
  }
\centering
\begin{tabular}{| c | c | c| c| }

 \hline
  setting  &  \multicolumn{1}{c|}{ ELDG} & \multicolumn{1}{c|}{ GELDG}
  \\ \hline
    mesh $\nu_{j+\frac12}$ &  $\sin(x_{j+\frac12})$   &  $\sin(x_{j+\frac12}) $  \\
  \hline
    adjoint $\alpha_j$     &          --         &    $\sin(x_j)$       \\
\hline
  \end{tabular}
\label{The related parameters settings for example2}
\end{table}
\begin{table}[!ht]
\footnotesize
\caption{1D transport equation with variable coefficients. $u_t + (\sin(x)u)_x =0$ with the initial condition $u(x,0) = 1$. $T=1$.
We use $CFL=0.3$ and $CFL=0.18$ for all $P^1$ (RK2) and $P^2$ (RK3) schemes, respectively.
  }
\centering
\begin{tabular}{| c | cc| cc|  }
\hline
 Mesh  &  {$L^1$ error} & Order &  {$L^1$ error} & Order \\
 \hline
  & \multicolumn{2}{c|}{$P^1$ GELDG} & \multicolumn{2}{c|}{$P^1$ GELDG}
  \\ \hline
    40 &      1.36E-03 &    --    &   1.36E-03   &   --\\
    80 &      3.57E-04 &     1.93 &   3.57E-04   &     1.93\\
   160 &      8.95E-05 &     1.99 &   8.95E-05   &     1.99\\
   320 &      2.31E-05 &     1.95 &   2.31E-05   &     1.95\\
\hline
  & \multicolumn{2}{c|}{$P^2$ ELDG} & \multicolumn{2}{c|}{$P^2$ GELDG}
  \\ \hline
    40 &    5.15E-05 & -- &     5.20E-05 & --\\
    80 &    6.33E-06 &     3.03 &     6.37E-06 &     3.03\\
   160 &    7.84E-07 &     3.01&     7.89E-07 &     3.01 \\
   320 &    9.60E-08 &     3.03&     9.71E-08 &     3.02 \\
\hline
\end{tabular}
\label{1dsin_example}
\end{table}

\begin{figure}[!ht]
\centering
\includegraphics[height=50mm]{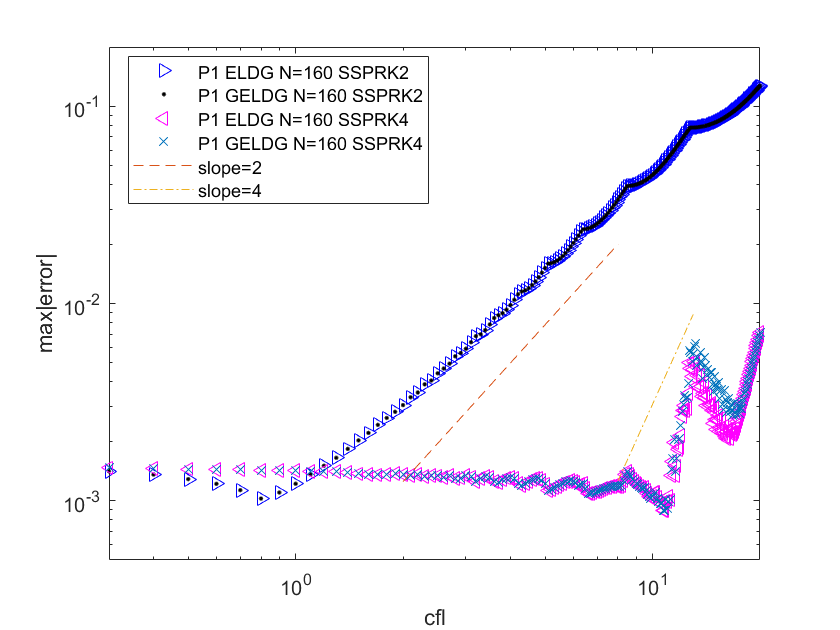}
\includegraphics[height=50mm]{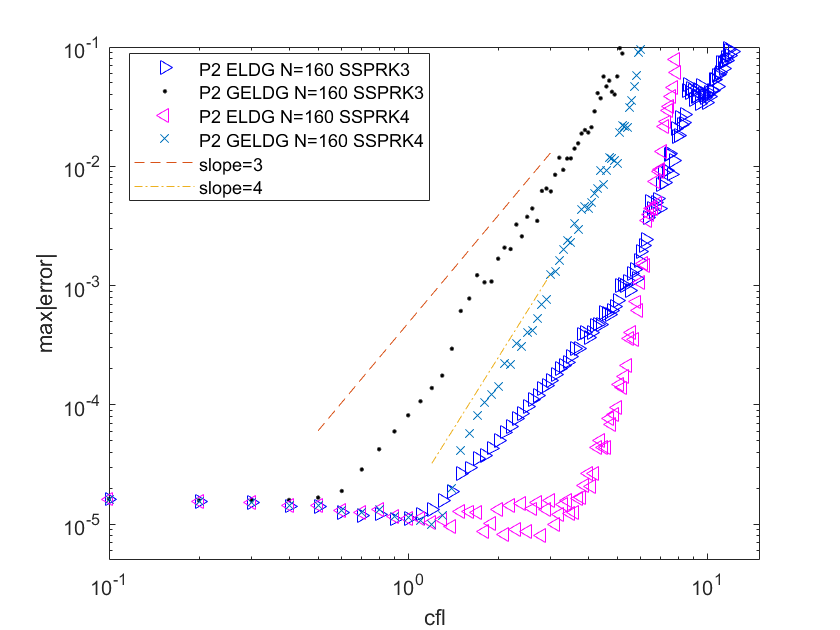}
\caption{The $L^\infty$ error versus $CFL$ of ELDG methods and GELDG methods for ~\eqref{eqn: 1D transport equation with variable coefficients} with the initial condition $u(x,0) = 1$. $T=1$.    $\Delta t=CFL\Delta x$.
}
\label{sin1d_time}
\end{figure}

\begin{table}[!ht]
\footnotesize
\caption{1D transport equation with variable coefficients. $u_t + (\sin(x)u)_x =0$ with the initial condition $u(x,0) = 1$ and inflow boundary condition. $T=1$.
We use $CFL=0.3$ and $CFL=0.18$ for all $P^1$ (RK2) and $P^2$ (RK3) schemes, respectively.
  }
\centering
\begin{tabular}{| c | cc| cc|  }
\hline
 Mesh  &  {$L^1$ error} & Order &  {$L^1$ error} & Order \\
 \hline
  & \multicolumn{4}{c|}{$P^1$ GELDG}\\
 \hline
  & \multicolumn{2}{c|}{without PP limiter} & \multicolumn{2}{c|}{with PP limiter}
  \\ \hline
    20 &      5.06E-03 &     --   &   4.97E-03   &     --\\
    40 &      1.36E-03 &     1.90 &   1.36E-03   &     1.88\\
    80 &      3.55E-04 &     1.94 &   3.55E-04   &     1.94\\
   160 &      8.90E-05 &     1.99 &   8.90E-05   &     1.99\\
   320 &      2.30E-05 &     1.95 &   2.30E-05   &     1.95\\
\hline
  & \multicolumn{4}{c|}{$P^2$ GELDG}\\
\hline
  & \multicolumn{2}{c|}{without PP limiter} & \multicolumn{2}{c|}{with PP limiter}
  \\ \hline
    20 &    4.16E-04 &     --   &     4.16E-04 &     --\\
    40 &    5.20E-05 &     3.00 &     5.20E-05 &     3.00\\
    80 &    6.38E-06 &     3.03 &     6.38E-06 &     3.03\\
   160 &    7.90E-07 &     3.01 &     7.90E-07 &     3.01 \\
   320 &    9.72E-08 &     3.02 &     9.72E-08 &     3.02 \\
\hline
\end{tabular}
\label{1dsin_example inflow boundary}
\end{table}
\begin{figure}[!ht]
\centering
\includegraphics[height=50mm]{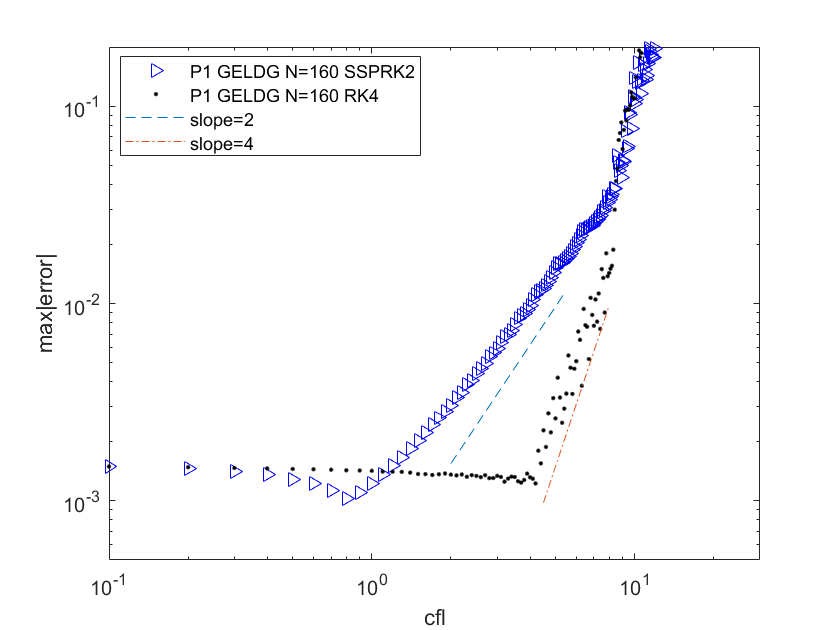}
\includegraphics[height=50mm]{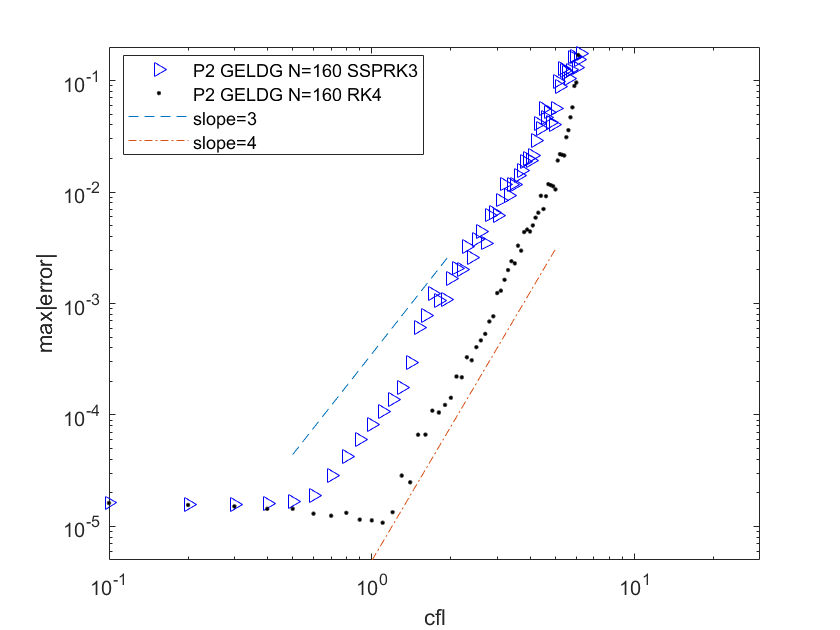}
\caption{The $L^\infty$ error versus $CFL$ of  GELDG method with PP limiter for ~\eqref{eqn: 1D transport equation with variable coefficients} with the initial condition $u(x,0) = 1$ and inflow boundary. $T=1$.    $\Delta t=CFL\Delta x$.
}
\label{sin1d_time inflow boundary}
\end{figure}

\end{exa}

\subsection{2D linear passive-transport problems}

\begin{exa}
(Swirling deformation flow). Consider
\begin{equation}\label{swirling deformation}
u_t - (\cos^2(\frac{x}{2})\sin(y)g(t)u)_x + (\sin(x)\cos^2(\frac{y}{2})g(t)u)_y =0, \ (x,y)\in[-\pi,\pi]^2
\end{equation}
where $g(t)=\cos(\frac{\pi t}{T})\pi$ and $T=1.5$. The initial condition is the following smooth cosine bells (with $C5$ smoothness),
\begin{equation}\label{cosine bell}
  u(x,y,0)=\left\{
  \begin{aligned}
    &r_0^b \cos^6(\frac{r^b}{2r_0^b}\pi), &if r^b<r_0^b,\\
    &0, \ &otherwise,
  \end{aligned}
  \right.
\end{equation}
where $r^b_0 = 0.3\pi$, and $r^b =\sqrt{(x-x_0^b)^2+(y-y_0^b)^2}$ denotes the distance between $(x, y)$ and the center of the cosine bell $(x_0^b, y_0^b) = (0.3\pi, 0)$. In fact, along the direction of the flow, the initial function becomes largely deformed at $t = T /2$, then goes back to its initial shape at $t = T$ as the flow reverses. If this problem is solved up to $T$, we call such a procedure one full evolution. We test accuracy for $Q^k$ EL DG, GEL DG methods without and with PP limiter with 4th RK and 4th splitting method for $k = 1, 2$ with $CFL = 2.5 $ up to one full evolution, and summarize results in Tables \ref{swirling deformation for cosine bell}. As expected, the $(k + 1)$th order convergence is observed for these methods. We plot the $L^\infty$ error versus $CFL$ of EL RK DG, GEL DG methods with $Q^1$ (left) and $Q^2$ (right) polynomial spaces for this case in Figure \ref{swirlingdeformation_consinebell}
which show that these three methods have similar stability for higher order discretization. We also plot the
$L^\infty$ error versus $CFL$ of GEL DG methods with PP limiter in Figure \ref{swirlingdeformation_consinebellpp}, in comparison to the one without limiter in the previous figure.

\begin{table}[!ht]\footnotesize
\caption{Swirling deformation flow. $Q^k$ EL DG and GEL DG methods without and with PPlimiter (k = 1, 2) for \eqref{swirling deformation} with
the smooth cosine bells \eqref{cosine bell} at T = 1.5. $CFL = 2.5$.}
\centering
\begin{tabular}{|c|c|c|c|c|c|c|c|c|c|}
\hline
\multirow{2}{*}{} & \multirow{2}{*}{N} & \multicolumn{4}{|c|}{$Q^1$} & \multicolumn{4}{|c|}{$Q^2$} \\
\cline{3-10}
& & $L^2$-error & order & $L^\infty$-error & order & $L^2$-error & order & $L^\infty$-error & order \\
\hline
\multirow{5}{*}{EL DG} & $20^2$ & 1.85E-02 & - & 2.78E-01 & - & 3.95E-03 & - & 5.74E-02 & -  \\
& $40^2$ & 4.14E-03 & 2.16 & 8.21E-02 & 1.76 & 1.76E-04 & 4.49 & 3.99E-03 & 3.85\\
& $80^2$ & 6.29E-04 & 2.72 & 1.39E-02 & 2.56 & 1.59E-05 & 3.47 & 3.57E-04 & 3.48\\
& $160^2$ & 9.05E-05 & 2.80 & 2.28E-03 & 2.61 & 2.12E-06 & 2.90 & 5.25E-05 & 2.77\\
& $320^2$ & 1.52E-05 & 2.57 & 4.19E-04 & 2.45 & 2.73E-07 & 2.96 & 6.82E-06 & 2.94\\
\hline
\multirow{5}{*}{GEL DG} & $20^2$ & 1.85E-02 & - & 2.78E-01 &- & 3.96E-03 & - & 5.78E-02 & -\\
& $40^2$ & 4.14E-03 & 2.16 & 8.21E-02 & 1.76 & 1.76E-04 & 4.49 & 4.01E-03 & 3.85 \\
& $80^2$ & 6.29E-04 & 2.72 & 1.39E-02 & 2.56 & 1.58E-05 & 3.47 & 3.54E-04 & 3.50 \\
& $160^2$ & 9.05E-05 & 2.80 & 2.28E-03 & 2.61 & 2.12E-06 & 2.90 & 5.25E-05 & 2.76 \\
& $320^2$ & 1.52E-05 & 2.57 & 4.19E-04 & 2.45 & 2.73E-07 & 2.96 & 6.82E-06 & 2.94\\
\hline
\multirow{5}{*}{GEL DG-PP} & $20^2$ & 2.07E-02 & - & 3.22E-01 & - & 3.61E-03 & - & 5.13E-02 & -\\
& $40^2$ & 4.16E-03 & 2.31 & 8.65E-02 & 1.89 & 2.96E-04 & 3.61 & 5.53E-03 & 3.21 \\
& $80^2$ & 6.34E-04 & 2.71 & 1.40E-02 & 2.63 & 1.95E-05 & 3.92 & 4.82E-04 & 3.52 \\
& $160^2$ & 9.01E-05 & 2.81 & 2.26E-03 & 2.63 & 2.15E-06 & 3.18 & 5.40E-05 & 3.16 \\
& $320^2$ & 1.51E-05 & 2.58 & 4.15E-04 & 2.44 & 2.74E-07 & 2.97 & 6.92E-06 & 2.96\\
\hline
\end{tabular}
\label{swirling deformation for cosine bell}
\end{table}

\begin{figure}[!ht]
\centering
\includegraphics[height=50mm]{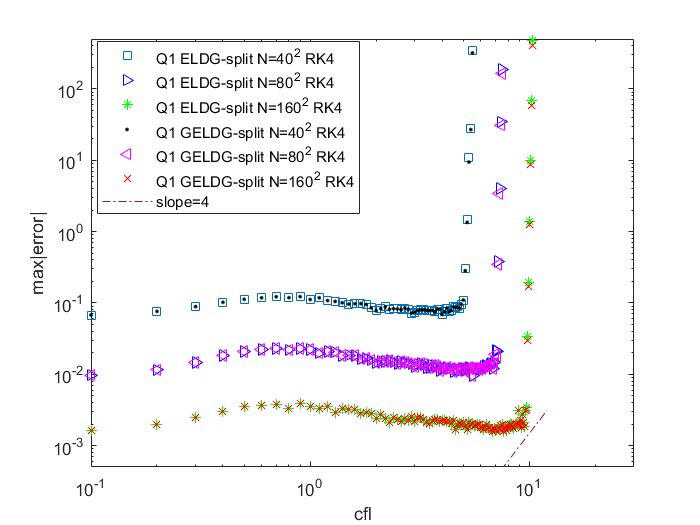}
\includegraphics[height=50mm]{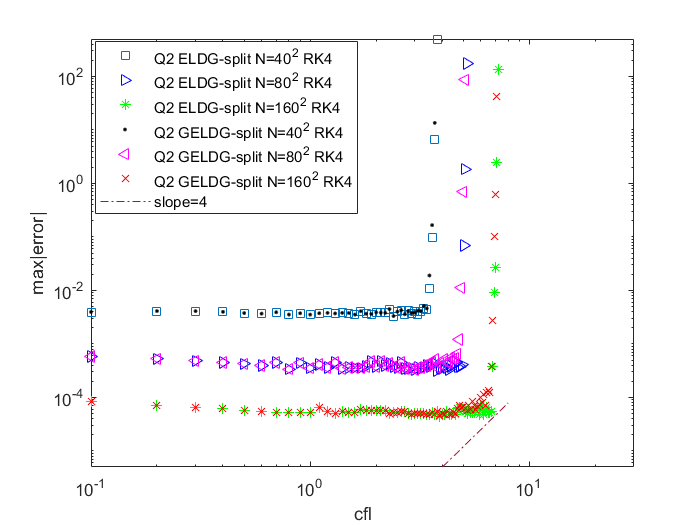}
\caption{The $L^\infty$ error versus $CFL$ of EL DG methods and GEL DG methods for \eqref{swirling deformation} with
the smooth cosine bells \eqref{cosine bell} at T = 1.5.
}
\label{swirlingdeformation_consinebell}

\end{figure}

\begin{figure}[!ht]
\centering
\includegraphics[height=50mm]{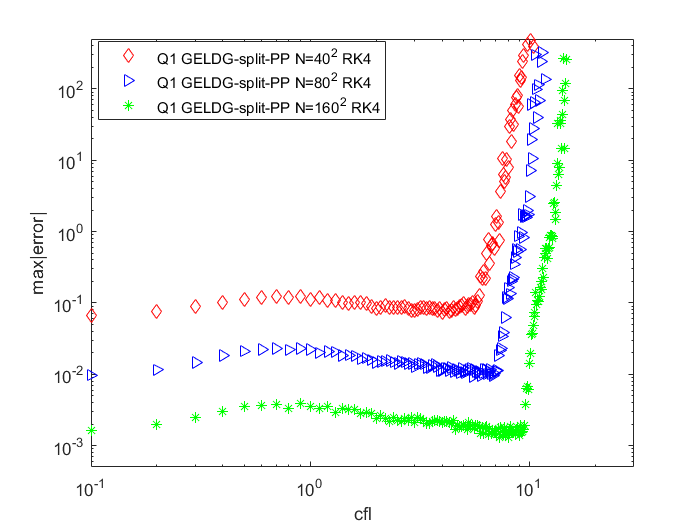}
\includegraphics[height=50mm]{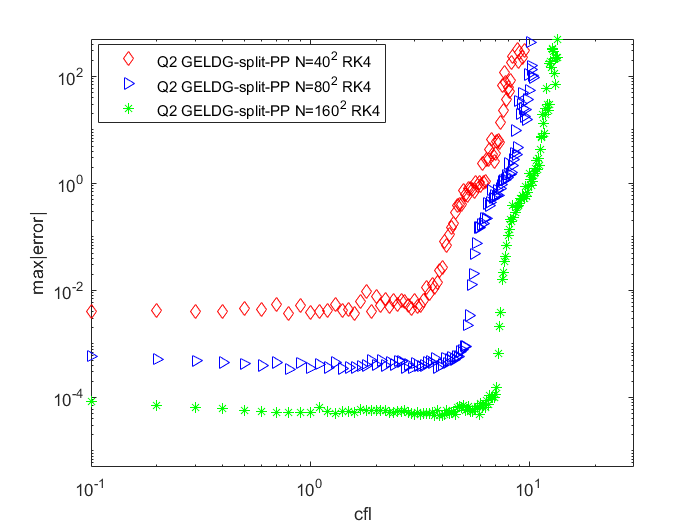}
\caption{The $L^\infty$ error versus $CFL$ of GEL DG with PP limiter methods for \eqref{swirling deformation} with
the smooth cosine bells \eqref{cosine bell} at T = 1.5.
}
\label{swirlingdeformation_consinebellpp}
\end{figure}

Lastly, we test two schemes on the swirling deformation flow \eqref{swirling deformation} with the following setting: the computational domain is $[- \pi, \pi]$ with the periodic boundary conditions and an initial condition plotted in Figure \ref{threeshape_initial_P2400}, which consists of a slotted disk, a cone as well as a smooth hump, similar to the one used in \cite{leveque1996high}. It is a challenging test for controlling oscillations around discontinuities. We adopt a simple TVB limiter with $M=15$ in \cite{cockburn1989TVB} for all schemes. We simulate this problem after one full revolutions and report the numerical solutions in Figure \ref{PPthreeshape3D_eldgandgeldg}. For better comparison, we plot 1D cuts and zoom in of the numerical solutions along with the exact solution to demonstrate the effectiveness of the PP limiter in Figure \ref{threeshape1D_eldgandgeldgPP}. It is found that oscillations are well controlled with the TVB limiter and are positivity preserving with the PP limiter. Solutions with larger CFL are observed to dissipate less than solutions from smaller CFL.

\begin{figure}[!ht]
\centering
\includegraphics[height=90mm]{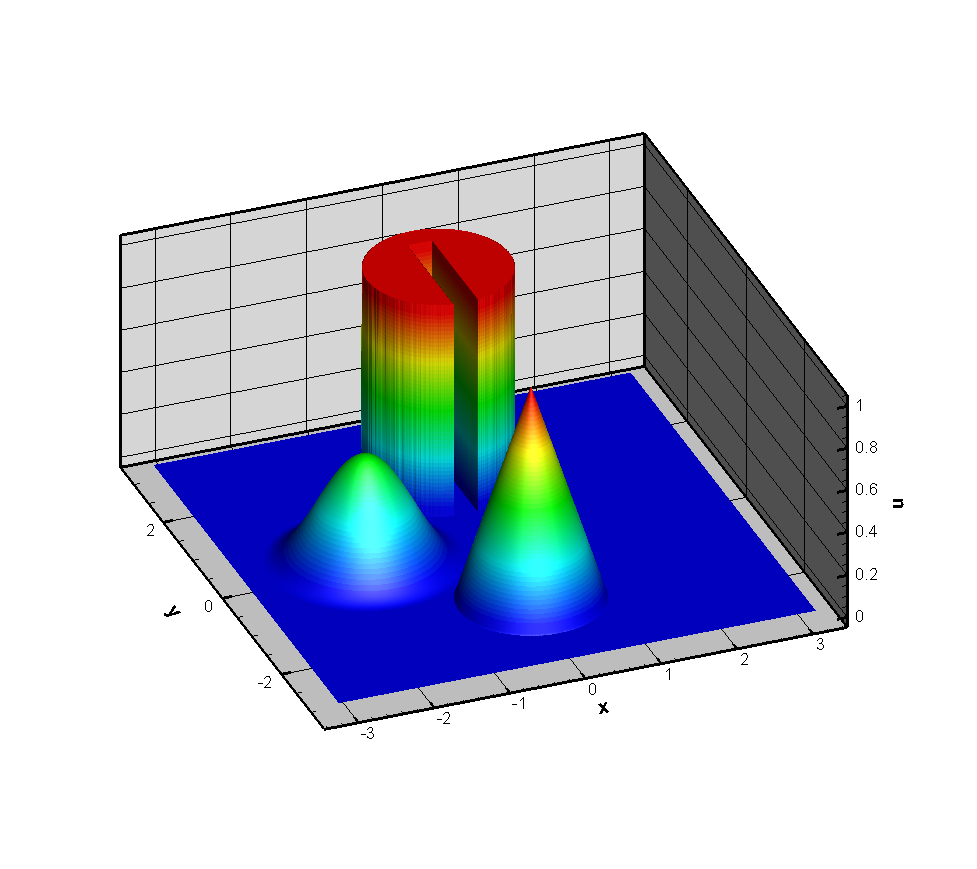}
\caption{Plots of the initial profile. The mesh of $400\times 400$ is used.
}
\label{threeshape_initial_P2400}
\end{figure}

\begin{figure}[!ht]
\centering
\includegraphics[height=60mm]{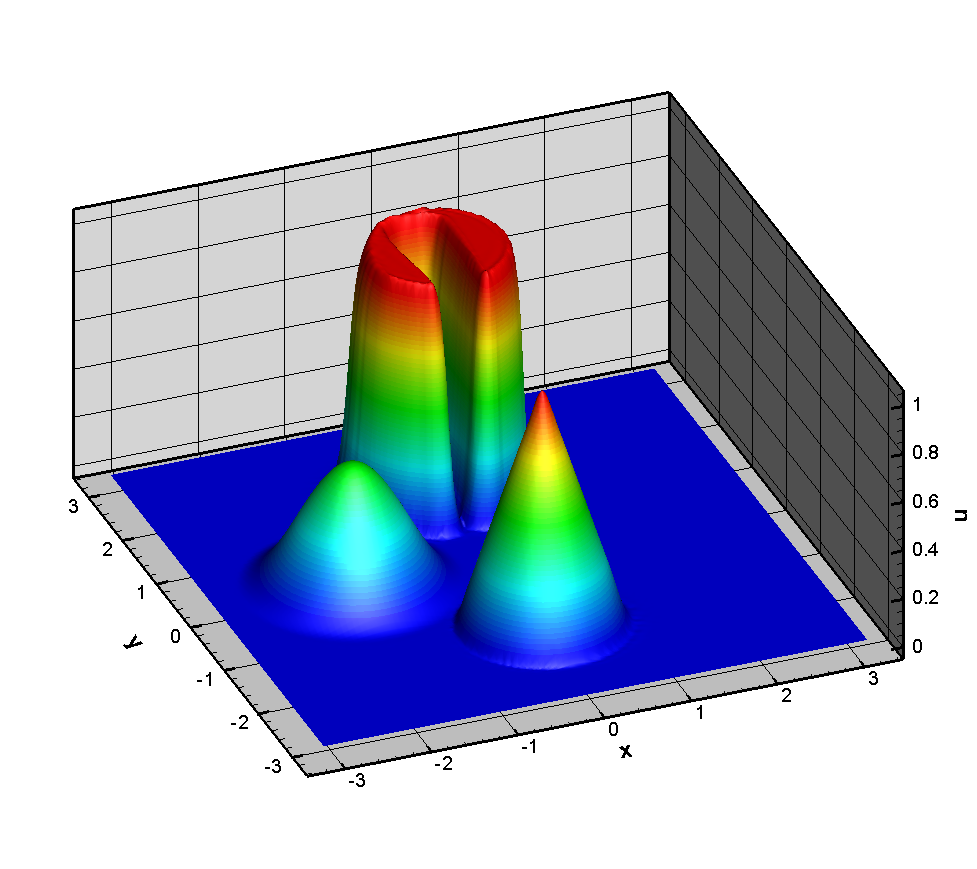}
\includegraphics[height=60mm]{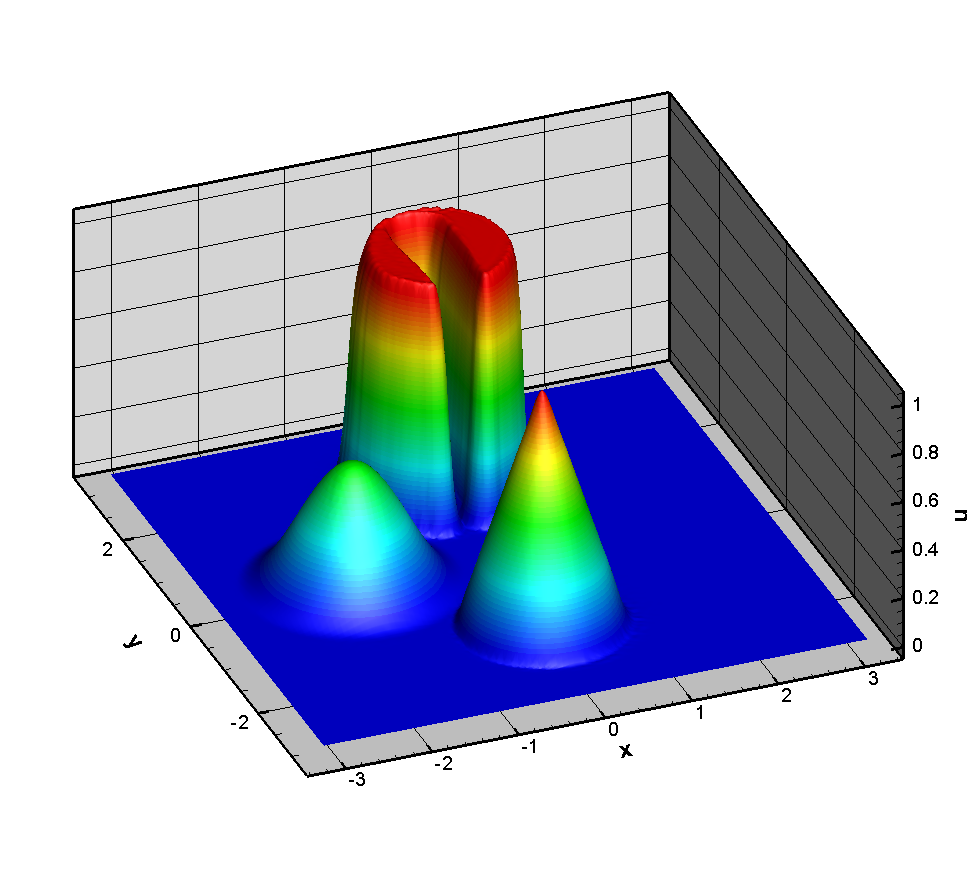}
\caption{Plots of the numerical solutions of GEL DG schemes with TVB and PP limiters for solving \eqref{swirling deformation} with initial data plotted in Fig. \ref{threeshape_initial_P2400}. The final integration time T is 1.5. The mesh of $100\times 100$ is used. Left: $Q^2$ GELDG-split+TVB+PPlimiter with $CFL = 2.2$. Right: $Q^2$ GELDG-split+TVB+PPlimiter with $CFL = 5.2$.
}
\label{PPthreeshape3D_eldgandgeldg}
\end{figure}



\begin{figure}[!ht]
\centering
\includegraphics[height=60mm]{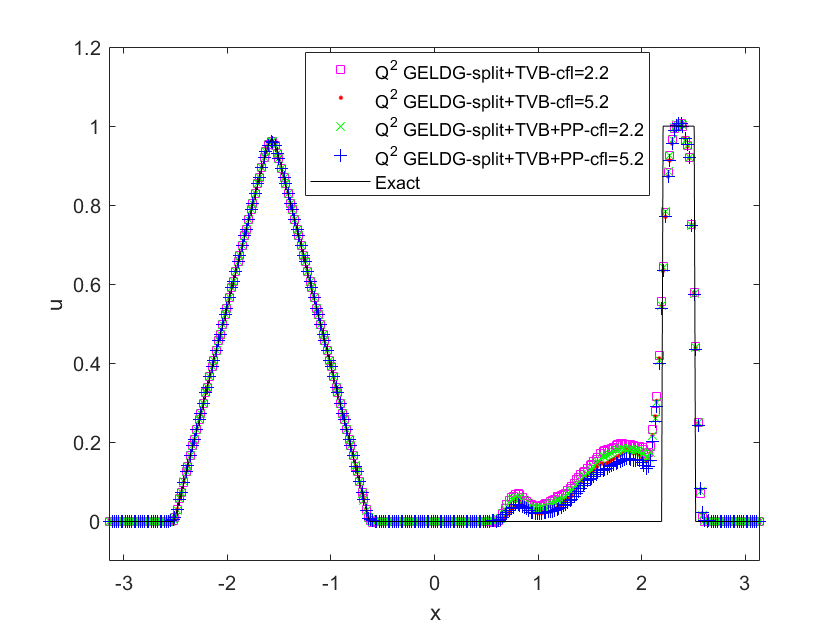}
\includegraphics[height=60mm]{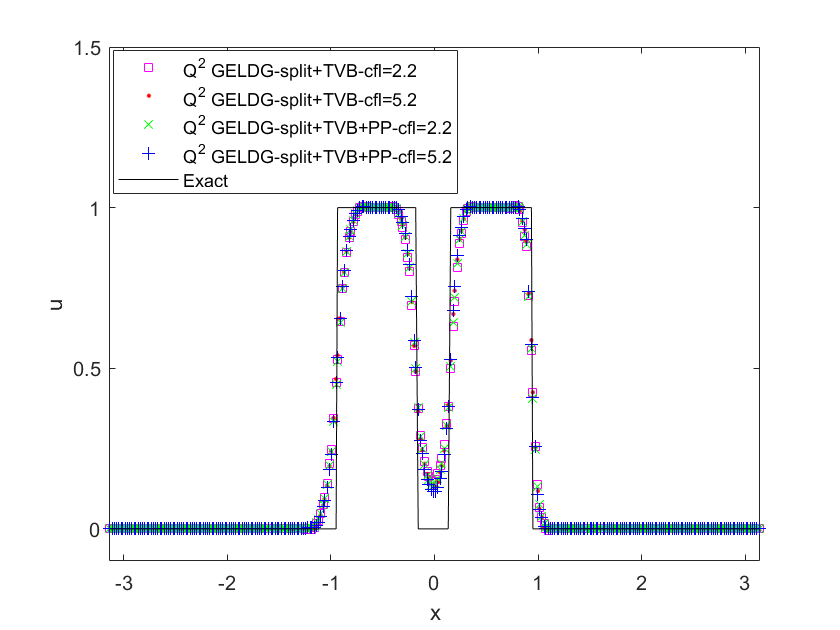}
\includegraphics[height=60mm]{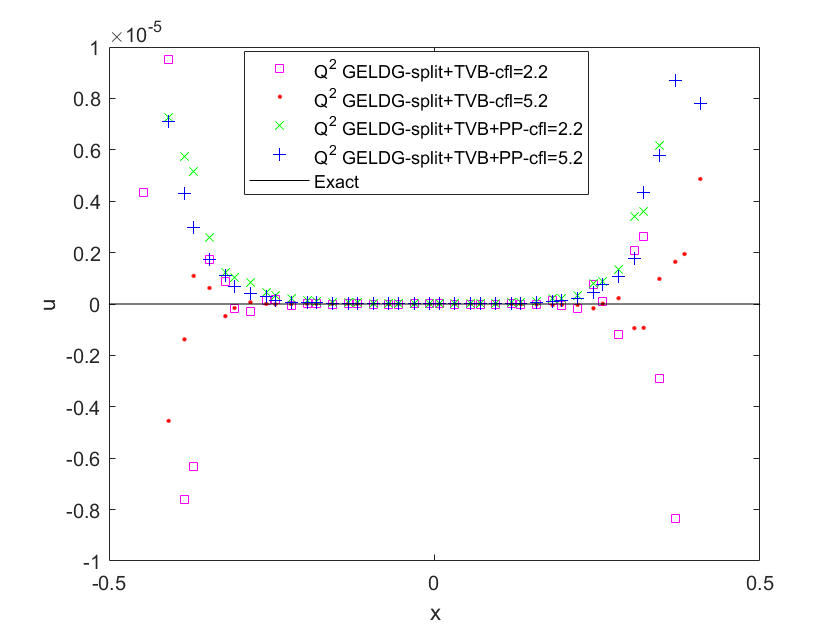}
\includegraphics[height=60mm]{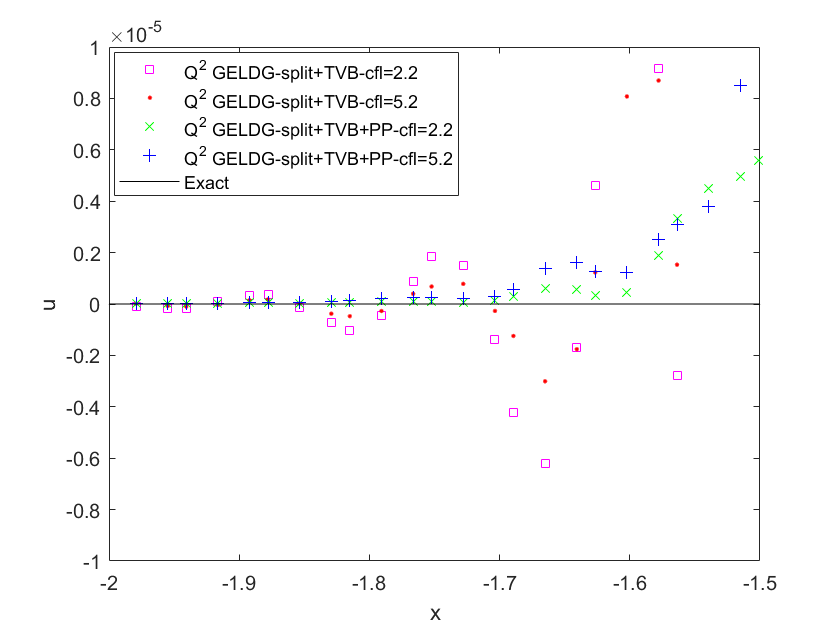}
\caption{Plots of 1D cuts and zoom in of the numerical solution for GEL DG methods with and without PP limiter for solving \eqref{swirling deformation} with initial data Fig. \ref{threeshape_initial_P2400}. The mesh of $100\times 100$ is used. Left: numerical solution at $x = 0 + \pi/100$. Right: numerical solution at $y = \pi/2 + \pi/100$.
}
\label{threeshape1D_eldgandgeldgPP}
\end{figure}


\end{exa}

\section{Conclusion} \label{section:conclusion}

In this paper, we develop a generalized Eulerian-Lagrangian (GEL) discontinuous Galerkin (DG) method for linear transport problems. Inflow boundary treatment is discussed. The method has the advantages in stability under large time stepping sizes, and in mass conservation, compactness and high order accuracy. Maximum principle preserving and positivity preserving limiters are proposed, leading to the discrete geometric conservation laws. These properties are numerically verified by extensive numerical tests for 1D and 2D linear transport equations. Future works include further theoretic development and developing schemes for linear system such as the wave equations.


\bibliographystyle{abbrv}
\bibliography{refer17}

\end{document}